\documentclass[preprint,12pt]{elsarticle}

\usepackage{geometry}
\usepackage{color}
\usepackage[utf8]{inputenc}
\usepackage[USenglish]{babel}
\usepackage[T1]{fontenc}

\usepackage{lmodern}
\usepackage{amsmath}
\usepackage{bm}
\usepackage{amssymb}
\usepackage{mathtools}
\usepackage{amsthm}
\usepackage{latexsym}
\usepackage{fixmath}
\usepackage{upgreek}

\usepackage{siunitx-v2}
\usepackage{xparse}
\usepackage{hyperref}
\usepackage{natbib}

\biboptions{sort&compress}

\usepackage{subcaption}
\usepackage[capitalise,]{cleveref}

\usepackage[font=small,labelfont=bf]{caption}
\usepackage{booktabs}
\usepackage{graphicx}

\makeatletter
\def\ps@pprintTitle{%
	\let\@oddhead\@empty
	\let\@evenhead\@empty
	\let\@oddfoot\@empty
	\let\@evenfoot\@oddfoot
}
\makeatother

\begin{document}

\begin{frontmatter}


\title{A new numerical \textit{mesoscopic scale} one-domain approach solver for free fluid/porous medium interaction}

\author[label1]{Costanza Aricò\corref{cor1}}

\ead{costanza.arico@unipa.it}
\cortext[cor1]{Corresponding Author}
\affiliation[label1]{organization={Department of Engineering, University of Palermo},
	addressline={viale delle Scienze},
	city={Palermo},
	postcode={90128},
	country={Italy}}

\author[label2]{Rainer Helmig}
\affiliation[label2]{organization={Institute for Modelling Hydraulic and Environmental Systems (IWS),
		             Department of Hydromechanics and Modelling of Hydrosystems, University of Stuttgart},
	addressline={Pfaffenwaldring 61},
	city={Stuttgart},
	postcode={D-70569},
	country={Germany}}

\author[label1]{Daniele Puleo}
\author[label2]{Martin Schneider}

\begin{abstract}
	A new numerical continuum \textit{one-domain} approach (ODA) solver is presented for the simulation of the transfer processes between a free fluid and a porous medium. The solver is developed in the \textit{mesoscopic} scale framework, where a continuous variation of the physical parameters of the porous medium (e.g., porosity and permeability) is assumed. The Navier-Stokes-Brinkman equations are solved along with the continuity equation, under the hypothesis of incompressible fluid. The porous medium is assumed to be fully saturated and can potentially be anisotropic. The domain is discretized with unstructured meshes allowing local refinements. A fractional time step procedure is applied, where one predictor and two corrector steps are solved within each time iteration. The predictor step is solved in the framework of a marching in space and time procedure, with some important numerical advantages. The two corrector steps require the solution of large linear systems, whose matrices are sparse, symmetric and positive definite, with $\mathcal{M}$-matrix property over Delaunay-meshes. A fast and efficient solution is obtained using a preconditioned conjugate gradient method. The discretization adopted for the two corrector steps can be regarded as a Two-Point-Flux-Approximation (TPFA) scheme, which, unlike the standard TPFA schemes, does not require the grid mesh to be $\mathbf{K}$-orthogonal, (with $\mathbf{K}$ the anisotropy tensor). As demonstrated with the provided test cases, the proposed scheme  correctly retains the anisotropy effects within the porous medium. Furthermore, it overcomes the restrictions of existing mesoscopic scale one-domain approachs proposed in the literature.
\end{abstract}


\begin{highlights}
\item Investigation of free fluid and porous media coupling transfer interface processes
\item New numerical solver in the framework of the mesoscopic one-domain approach
\item Incompressible fluid within and around an isotropic or anisotropic porous medium
\item Numerical solution of the Navier-Stokes-Brinkman equations by a fractional time step procedure
\item Some real-world applications are presented
\end{highlights}


\end{frontmatter}


\section{Introduction}	\label{intro}  
Momentum transfer at the interface between free fluid and porous media is of significance for various applications. Indeed, interface transport processes are involved in different industrial, environmental and biological/biomedical applications, as for example passive control devices using porous coating, heat exchangers, fuel cells, filtration and drying processes, groundwater pollution, flows in fractured media, geothermal systems, flows in biological tissues and related medical drugs transport problems. The study of fluid-porous interface momentum transfer is also crucial for the development of mathematical and numerical models involving additional transfer processes, e.g., passive solute or heat and pollutant transport. \par 
The study of the fluid–porous interface transfer can be performed at different scales \cite{1}. At the microscopic pore scale, the flow in the free fluid region and in the void spaces of the  porous medium is governed by the classical (Navier)-Stokes equations along with boundary conditions at the interface between the fluid phase and the solid phase within the permeable region (e.g., no-slip velocity condition). Such a pore scale approach has two limitations, 1) only small-scale problems can be simulated due to the high computational effort, caused by the discretization of the microscopic void spaces, wich requires a significant number of mesh elements, and 2) often, detailed knowledge of the pore geometry in the entire domain is not known, even with advanced image acquisition technologies. This is the reason why descriptions at the mesoscopic and macroscopic scales are usually introduced according to two popular approaches, both derived from the volume averaging of pore-scale governing equations. \par
In the \textit{two-domain} approach (TDA), at the macroscopic scale, the bulk fluid and porous regions are separated by a sharp interface. Two different sets of governing equations are applied in each of the bulk regions, namely the Darcy equation in the porous domain and the (Navier)-Stokes equations in the free fluid domain. Due to the different character of the corresponding partial differential equations, specific boundary conditions have to be imposed at the common interface that guarantee conservation of the fluxes together with appropriate slip conditions for tangential component of the free fluid velocity \cite{2}. In the pioneering work of Beavers and Joseph \cite{3}, such a slip boundary condition was experimentally derived for parallel flow conditions when coupling Stokes with Darcy flow. Ochoa-Tapia and Whitaker \cite{4, 5} coupled the Stokes and Brinkman equations at the interface, assuming continuous tangential velocity but discontinuous tangential shear stress. Other interface boundary conditions have been obtained by using homogenization methods, see for example \cite{6, 7, 8} and references therein. \par
\begin{figure} [h]
	\centering
	\includegraphics[width=\textwidth]{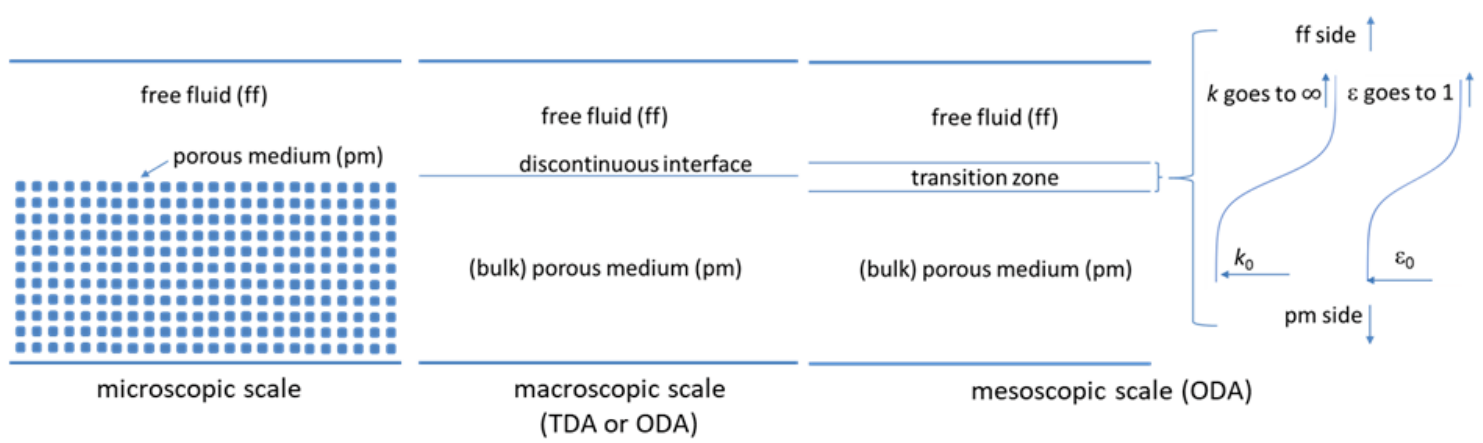}
	\caption{Pore scale reference configuration (left) and different REV-scale approaches. Sharp interface approaches are depicted in the middle picture (TDA or macroscopic ODA), whereas the ODA that uses a transition zone is shown in the right picture. There, some possible transition functions are plotted, where $k_0$ and $\epsilon_0$ denote the porosity and permeability values of the bulk porous medium}
	\label{fig:Figure 1}		
\end{figure}

In the continuum \textit{one-domain} approach (ODA), a ``fictitious'' equivalent single medium replaces the fluid and solid phases, and one set of governing equations, valid everywhere in the domain, is used to model the transfer processes. At the mesoscopic scale, the transition from the bulk fluid to the bulk porous region is modelled using a transition zone (or transition layer, TL) located in between these two regions \cite{1, 9, 10}. Within this transition layer the change of effective macroscopic properties of the permeable medium, such as porosity $\epsilon$ and permeability \textit{K}, is modelled with appropriate continuous transition functions. The set of governing equations is often denoted as (Navier-)Stokes-Brinkman, (Navier-)Stokes-Darcy, Darcy-Brinkmann or Brinkman equations \cite{9, 10}. These are derived by averaging the governing pore-scale equations over a Representative Elementary Volume (REV) \cite {11, 12}. The REV characteristic size is much smaller than that of the investigated domain, but much larger than the characteristic pore-scale size. In \cite{13} the authors solve the transfer momentum problem in a three-layer 1D channel, including two external bulk homogeneous porous and fluid regions, separated by a heterogeneous transition zone, with variable $\epsilon$ and \textit{K}. They give analytical velocity expressions in the three zones. In \cite{14} the authors present an ODA based on the Stokes-Brinkman model. By using the method of matched asymptotic expansion, the asymptotic solution for vanishing transition layer thickness is investigated. In \cite{15, 16}, an ODA is presented for 1D problems, where a macroscopic momentum equation, with Darcy form, applicable everywhere in the system, is solved along a homogenization closure problem, to obtain the transition layer permeability profile. The porous medium is assumed to be periodic and periodicity conditions are imposed for the closure problem. The ODA model is derived under the assumptions of constant pressure gradient and a given convective fluid velocity in the inertial term of the momentum equation. 
Another macroscopic ODA approach uses penalization, such that the so-called \textit{penalized} Navier-Stokes equations are applied, with an extra penalizing Darcy term in the momentum equation, which accounts for the drag force of the solid particles of the porous medium over the fluid \cite{17, 18, 19, 20, 21, 22}. The Darcy term is a function of the porosity and inverse of the permeability. This term is applied only within the porous region, while in the clear fluid region it vanishes, such that the classical (Navier)-Stokes equations are solved there. This implies that continuous transition of  porosity and permeability is assumed close to the free fluid-porous medium interface, but instead a discontinuous change is considered. These discontinuities induce an interfacial stress jump, and an additional stress arises within the porous medium due to the Darcy term in the governing momentum  equations. Such penalized approaches have been widely applied in the literature since they adopt well-consolidated numerical procedures for the classical (Navier)-Stokes equations. These different approaches are schematically depicted in \cref{fig:Figure 1} together with the reference configuration on the pore scale. \par 
Each of these approaches has advantages and drawbacks. If on the one hand the sets of governing equations in the TDA are well-known and consolidated numerical tools can be applied, a good match of the results near the interface is obtained when choosing appropriate interface coupling conditions \cite{9}. Many applications present a gradual variation of the macroscopic properties (porosity, permeability, etc.) of the porous medium, without any abrupt change between the two bulk fluids and permeable regions \cite{9}, so the ODA should be more suitable. The principal limitation of the ODA is related to the difficulty in predicting the spatial variation of $\epsilon$ and \textit{K} in the transition layer. \par 
In literature, most of the mesoscopic ODA models are analytically solved under specific geometrical conditions and by assuming simplified boundary conditions, e.g., 1D flow, periodicity of the flow and the porous medium, steady-state flow or Stokes flow regime (low Reynolds number), see for example \cite{13, 14, 15, 16}. Several numerical procedures have been proposed in the literature for the macroscopic ODA and TDA. In \cite{23} and references therein, finite element methods, mixed methods and discontinuous Galerkin methods, as well as their combinations, are discussed. All these procedures involve solving different models in each of the bulk regions, and couple them using suitable boundary conditions. Numerical challenges when solving Navier-Stokes-Brinkman equations are discussed in \cite{23}. These are mainly related to the coupling of Galerkin approximations of both Stokes and Darcy problems with mixed Finite Element formulations. Indeed, combining  Raviart–Thomas Finite Element velocity spaces \cite{24} with piecewise constant or linear pressure fields can satisfy the \textit{inf-sup} conditions \cite{23}. In \cite{25}, such discretization allows us to obtain the correct solution for the Darcy equation, but is not suitable for the Stokes problems. An alternative is to invoke the Darcy's law in the mass conservation equation, leading to an elliptic pressure Poisson problem, which can be easily approximated by Galerkin techniques. Unfortunately, this technique leads to a loss of accuracy for the velocity solution, as well as to a weak enforcement of the mass conservation equation \cite{25}. Stabilised Finite Element methods (e.g.,  Galerkin/Least-Squares methods, Streamline-Upwind/Petrov–Galerkin or Pressure-Stabilising/Petrov–Galerkin methods, Pressure Gradient Projection methods or Variational Multi-Scale methods) have been successfully applied either for Darcy or Navier-Stokes equations (\cite{22} and the references therein). \par 
In the present paper, we propose a new numerical ODA solver at the mesoscopic level, where we do not restrict the porous medium to be isotropic, but also consider anisotropy. We solve the continuity and the Navier-Stokes-Brinkman equations, applying a fractional time step procedure, where a prediction and two correction problems are sequentially solved within each time iteration. The computational domain is discretized using unstructured meshes, allowing mesh refinement at the transition layer. The solution of the prediction problem is performed by a Marching in Space and Time (MAST) procedure. This is a Finite Volume algorithm, recently presented for the solution of shallow water and groundwater problems, as well as Navier-Stokes flow applications (see for example \cite{26, 27, 28, 29, 30} and references therein). This scheme has some important features, which are further discussed in the following Sections (see also \cite{29, 30} and references therein). The two correction problems involve a fast and efficient solution of large linear systems since the associated matrices are sparse, symmetric, positive definite and diagonally dominant. As discussed in the following Sections, the algorithm proposed for the discretization of the correction problems can be regarded as a Two-Point-Flux Approximation (TPFA) scheme, which retains the anisotropic properties of the porous medium, but, unlike the standard TPFA scheme, it does not require the computational grid to be aligned with the principal anisotropy directions. The proposed method is \textit{strongly} conservative, in the sense that the velocity solution is divergence free pointwise inside each mesh cell, and local and global mass balance are always guaranteed. The method is suitable for simulation of multi-dimensional unsteady flow problems. \par    
The paper is organized as follows. The governing equations and the characteristics of the discretizing mesh are presented in \cref{sect_gov_eqq}, in \cref{numerical_model} we provide the algorithmic details of the new ODA solver, and in \cref{num_tests} some numerical applications, including the analysis of the convergence order and the computational costs, as well as some ``real-world'' applications, are presented. 

\section{Governing Equations} \label{sect_gov_eqq}
We assume a Newtonian incompressible fluid with density $\rho_l$ inside and around a saturated porous medium, which is assumed to be rigid, with solid particles fixed in space. At the mesoscopic scale, the fluid and solid phases are described as a single continuous medium, derived by averaging the micro-scale Navier-Stokes equations over a REV \cite{11, 12}. This yields the following governing equations \cite{1, 10}
\begin{subequations}  \label{eq:governing_Eqq}
	\begin{gather}
		\nabla\cdot\left(\rho_l \mathbf{u}\right)=0, \label{eq:continuity} \\
		\rho_l \frac{\partial \mathbf{u}}{\partial t}+\rho_l \mathbf{u}\cdot\nabla\left(\frac{\mathbf{u}}{\epsilon}\right) =-\epsilon \left(\nabla p^l - \rho_l \mathbf{g}\right)+ \mu \nabla\cdot\left(\nabla\mathbf{u}+\left(\nabla\mathbf{u}\right)^{T}\right) - \mu \epsilon \mathbf{\mathfrak{K}} \mathbf{u}, \label{eq:momentum}
	\end{gather}
\end{subequations}
\noindent where $t$ is time, $\mathbf{x}$ is the spatial coordinate vector, $\mathbf{u}$ is the surface average fluid velocity, with $u$ and $v$ its $x$ and $y$ components, $p^l$ is the intrinsic averaged fluid pressure \cite{10}, $\mathbf{g}$ is the gravitational acceleration, downward oriented, with $g$ the absolute value of its vertical component, $\mu$ is the dynamic fluid viscosity, $\epsilon$ is the porosity of the porous medium, and $\mathbf{\mathfrak{K}}$ is the inverse of the  permeability tensor of the porous medium, $\mathbf{K}$, symmetric and positive definite. The last term on the r.h.s. of \cref{eq:momentum} represents a drag force due to the microscopic momentum exchange of the fluid with the solid particles of the permeable matrix. According to \cite{4, 10}, it is related to $\mu$, to the relative velocity between the fluid and the solid grains and to the permeability of the porous  medium. \par 
Dividing \cref{eq:continuity,eq:momentum} by $\rho_l$ and setting $\Psi = \frac{p^l}{\rho_l} - gy$, we obtain 
\begin{subequations}  \label{eq:governing_Eqq2}
	\begin{gather}
		\nabla\cdot\left(\mathbf{u}\right)=0 , \label{eq:continuity2} \\
		\frac{\partial \mathbf{u}}{\partial t}+ \mathbf{u}\cdot\nabla\left(\frac{\mathbf{u}}{\epsilon}\right) =-\epsilon \left(\nabla \Psi \right)+ \nu \nabla\cdot\left(\nabla\mathbf{u}+\left(\nabla\mathbf{u}\right)^{T}\right) - \nu \epsilon \mathbf{\mathfrak{K}}  \mathbf{u}, \label{eq:momentum2}
	\end{gather}
\end{subequations}
\noindent with the kinematic fluid viscosity $\nu=\mu/\rho_l$. We solve the system in \cref{eq:continuity,eq:momentum} for the unknowns $\mathbf{u}$ and $\Psi$, in the computational domain $\Omega $, and let $\Gamma$ be its boundary. Three types of boundary conditions (BCs) can be assigned over $\Gamma=\Gamma_e\cup\Gamma_n\cup\Gamma_f$. $\Gamma_e$ is the portion where we assign \textit{essential} BCs (i.e., Dirichlet BCs for the velocity), $\Gamma_n$ the portion where we assign \textit{natural} BCs, (i.e., boundary stress vectors), and $\Gamma_f$ the portion where we assign \textit{free-slip} BCs, a combination of the previous ones. In \cref{eq:IBCs} we formulate the boundary and initial conditions (ICs) needed for the solution of system in \cref{eq:continuity,eq:momentum} to be well-posed.
\begin{subequations} \label{eq:IBCs}
	\begin{gather}
		\mathbf{u\left(\mathbf{x}\right)}=\mathbf{u}_b \qquad \forall \mathbf{x}\in\Gamma_e, \qquad t\geq0, \\
		\boldsymbol{\sigma\left(\mathbf{x}\right)}=\left(\Psi - 2\nu \frac{\partial u}{\partial n}  \right)\mathbf{n} - \nu \left(\frac{\partial u_n}{\partial s} + \frac{\partial u_s}{\partial n}\right) \mathbf{s} \qquad \forall \mathbf{x}\in\Gamma_{e} , \qquad t\geq0, \label{eq:visc}\\
		u_n=0 \qquad \boldsymbol{\tau}\left(\mathbf{x}\right) = 0 \qquad \forall \mathbf{x}\in\Gamma_{f} , \qquad t\geq0, \\
		\mathbf{u}=\mathbf{u}_0 \qquad \textrm{with} \qquad \nabla \cdot \mathbf{u}_0=0 \qquad \Psi=\Psi_0 \qquad \forall  \mathbf{x}\in\Omega , \qquad t= 0,
	\end{gather}
\end{subequations} 
\noindent where $\mathbf{u}_b$ and $\boldsymbol{\sigma}$ are the velocity and stress vectors imposed at the boundary, respectively, $u_n$ and $u_s$ are the components of $\mathbf{u}$ along the direction $\mathbf{n}$ and $\mathbf{s}$, normal (outward oriented) and tangent to the boundary, respectively, $\boldsymbol{\tau}$ is the stress vector along direction $\mathbf{s}$, and sub-index 0 marks the initial values of $\mathbf{u}$ and $\Psi$ in $\Omega$. In the rest of the paper, the viscous terms in \cref{eq:visc} are neglected. \par	
For the fractional time step procedure presented in the next Section, it is beneficial to re-write the system in \cref{eq:momentum2} in vector-matrix form as
\begin{equation} \label{eq:sys1}
	\frac{\partial \mathbf{U}}{\partial t} + \nabla\cdot\mathbf{F}\left(\mathbf{U}\right) + \mathbf{D}\left(\mathbf{U}\right)=\mathbf{B}\left(\mathbf{U}\right)+\nabla\cdot\mathbf{E}\left(\textbf{U}\right) ,
\end{equation}
\noindent where
\begin{subequations} \label{eq:sys2} %
	\begin{gather}
		\begin{gathered}
			\mathbf{U}=\begin{pmatrix}
				u\\v
			\end{pmatrix}\quad
			\mathbf{F\left(\mathbf{U}\right)}=\begin{pmatrix}
				\mathbf{F}_1 \quad \mathbf{F}_2
			\end{pmatrix}\quad
			\mathbf{E\left(\mathbf{U}\right)}=\begin{pmatrix}
				\mathbf{E}_1 \quad \mathbf{E}_2
			\end{pmatrix} \\
			\mathbf{B\left(\mathbf{U}\right)}=-\begin{pmatrix}
				\epsilon \nabla_{x}\Psi\\\epsilon \nabla_{y}\Psi
			\end{pmatrix} \quad
			\mathbf{D\left(\mathbf{U}\right)}= \epsilon\nu \mathbf{\mathfrak{K}} \begin{pmatrix}
				u\\v
			\end{pmatrix}
		\end{gathered},\\
		\mathbf{F}_1=\begin{pmatrix}
			u\frac{u}{\epsilon}\\u\frac{v}{\epsilon}
		\end{pmatrix}\quad
		\mathbf{F}_2=\begin{pmatrix}
			v\frac{u}{\epsilon}\\v\frac{v}{\epsilon}
		\end{pmatrix}\quad
		\mathbf{E}_1=\nu\begin{pmatrix}
			\frac{\partial u}{\partial x}\\ \frac{\partial v}{\partial x}
		\end{pmatrix}\quad
		\mathbf{E}_2=\nu\begin{pmatrix}
			\frac{\partial u}{\partial y}\\ \frac{\partial v}{\partial y}
		\end{pmatrix}.
	\end{gather}
\end{subequations}

\section{Numerical algorithm} \label{numerical_model}   

In \cref{FTS_procedure} we present a general overview of the fractional time step procedure applied to solve system in  \crefrange{eq:governing_Eqq}{eq:sys2}, while we refer to \cref{algo_details} those readers interested in the numerical details of the algorithm steps. 

\subsection{Fractional time step procedure} \label{FTS_procedure}  	
System in \cref{eq:governing_Eqq} is solved by applying a fractional time step procedure, where one predictor and two corrector problems are solved sequentially. This can be done by using the following splitting in \cref{eq:sys1,eq:sys2}
\begin{equation} \label{eq:fts}
	\begin{gathered}
		\mathbf{B}\left(\mathbf{U}\right) = \mathbf{B}\left(\mathbf{U}\right)^k+\left(\mathbf{B}\left(\mathbf{U}\right)-\mathbf{B}\left(\mathbf{U}^k\right)\right)\\  
		\mathbf{E}\left(\mathbf{U}\right) = \mathbf{E}\left(\mathbf{U}^{k-1/3}\right)+ \left(\mathbf{E}\left(\mathbf{U}\right)-\mathbf{E}\left(\mathbf{U}^{k-1/3}\right)\right)
	\end{gathered},
\end{equation}
\noindent such that system in \cref{eq:sys1} splits into
\begin{subequations} \label{FTS2}
	\begin{gather}
		\frac{\partial \mathbf{U}}{\partial t} + \nabla\cdot\mathbf{F}\left(\mathbf{U}\right)+ \mathbf{D}\left(\mathbf{U}\right)=\mathbf{B}\left(\mathbf{U}^k\right)+\nabla\cdot\mathbf{E}\left(\mathbf{U}^{k-1/3}\right) , \label{eq:Ps}\\
		\frac{\partial \mathbf{U}}{\partial t} + \mathbf{D}\left(\mathbf{U}\right) - \mathbf{D}\left(\mathbf{U}^{k+1/3}\right)= \nabla\cdot\mathbf{E}\left(\mathbf{U}\right) - \nabla\cdot\mathbf{E}\left(\mathbf{U}^{k-1/3}\right),  \label{eq:Cs1}\\
		\frac{\partial \mathbf{U}}{\partial t} + \mathbf{D}\left(\mathbf{U}\right) - \mathbf{D}\left(\mathbf{U}^{k+2/3} \right)= \mathbf{B}\left(\mathbf{U}\right) - \mathbf{B}\left(\mathbf{U}^k\right) , \label{eq:Cs2}
	\end{gather}
\end{subequations}
\noindent where \cref{eq:Ps} is the predictor problem (PP) and \cref{eq:Cs1,eq:Cs2} are the 1$^{st}$ and 2$^{nd}$ corrector problems (CP1 and CP2), respectively. In the following Sections, the symbols $t^k$, $t^{k+1/3}$, $t^{k+2/3}$ and $t^{k+1}$ mark the beginning of the time step, as well as the end of PP, CP1 and CP2, respectively. The symbol $t^{k-1/3}$ marks the end of the CP2 of the previous time step. According to a functional analysis, the system in \cref{eq:Ps} is a convective problem, while systems in \cref{eq:Cs1,eq:Cs2} are diffusive problems (see \cite{28, 29, 30} and literature therein). The time integral forms of \cref{eq:Ps,eq:Cs1,eq:Cs2} are
\begin{subequations} \label{int_form_FTS}
	\begin{gather}
		\begin{gathered}
			\mathbf{U}^{k+1/3}-\mathbf{U}^k+\nabla\cdot\int_0^{\Delta t}\mathbf{F}\left(\mathbf{U}\right)dt + \mathbf{D}\left(\mathbf{U}^{k+1/3}\right)\Delta t  =  \mathbf{B}\left(\mathbf{U}^k\right)\Delta t + \\ 
			\nabla \cdot \mathbf{E}\left(\mathbf{U}^{k-1/3}\right)\Delta t  \label{eq:int_f_1}
		\end{gathered}, \\
		\begin{gathered}
			\mathbf{U}^{k+2/3}-\mathbf{U}^{k+1/3}+ 
			\mathbf{D}\left(\mathbf{U}^{k+2/3}\right)\Delta t - \mathbf{D}\left(\mathbf{U}^{k+1/3}\right)\Delta t= \\\nabla \cdot \mathbf{E}\left(\mathbf{U}^{k+2/3}\right)\Delta t - \nabla \cdot \mathbf{E}\left(\mathbf{U}^{k-1/3}\right)\Delta t  \label{eq:int_f_2}
		\end{gathered},  \\
		\begin{gathered}
			\mathbf{U}^{k+1}-\mathbf{U}^{k+2/3} + \mathbf{D}\left(\mathbf{U}^{k+1}\right)\Delta t - \mathbf{D}\left(\mathbf{U}^{k+2/3}\right)\Delta t= \\ \mathbf{B}\left(\mathbf{U}^{k+1}\right)\Delta t - \mathbf{B}\left(\mathbf{U}^k\right)  \Delta t \label{eq:int_f_3} \\	      
		\end{gathered}	,		
	\end{gather}
\end{subequations}
\noindent and the sum of \crefrange{eq:int_f_1}{eq:int_f_3} gives the integral form of the original system in \cref{eq:sys1}. \par
The time discretization form of \crefrange{eq:Ps}{eq:Cs2} are
\begin{subequations} \label{discr_form_FTS}
	\begin{gather}
		\begin{gathered}
			\frac{\left(\mathbf{U}^{k+1/3}-\mathbf{U}^k\right)}{\Delta t} + \nabla\cdot \mathbf{F}\left(\mathbf{U}^{k+1/3}\right) + \mathbf{D}\left(\mathbf{U}^{k+1/3}\right) =  \mathbf{B}\left(\mathbf{U}^k\right) + \nabla \cdot \mathbf{E}\left(\mathbf{U}^{k-1/3}\right)   \label{eq:discr_f_1}
		\end{gathered}, \\
		\begin{gathered}
			\frac{\left(\mathbf{U}^{k+2/3}-\mathbf{U}^{k+1/3}\right)}{\Delta t}+ 
			\mathbf{D}\left(\mathbf{U}^{k+2/3}\right) - \mathbf{D}\left(\mathbf{U}^{k+1/3}\right) = \\\nabla \cdot \mathbf{E}\left(\mathbf{U}^{k+2/3}\right) - \nabla \cdot \mathbf{E}\left(\mathbf{U}^{k-1/3}\right)  \label{eq:discr_f_2}
		\end{gathered},  \\
		\begin{gathered}
			\frac{\left(\mathbf{U}^{k+1}-\mathbf{U}^{k+2/3}\right)}{\Delta t}+ \mathbf{D}\left(\mathbf{U}^{k+1}\right) - \mathbf{D}\left(\mathbf{U}^{k+2/3}\right)= \\ \mathbf{B}\left(\mathbf{U}^{k+1}\right) - \mathbf{B}\left(\mathbf{U}^k\right) \label{eq:discr_f_3} \\	      
		\end{gathered},			
	\end{gather}
\end{subequations}
We write in \cref{eq:discr_f_1} $\mathbf{D}\left(\mathbf{U}^{k+1/3}\right) = \left(\mathbf{D}\left(\mathbf{U}^{k+1/3}\right) - \mathbf{D}\left(\mathbf{U}^k\right)\right) + \mathbf{D}\left(\mathbf{U}^k\right)$, inserting in system in \cref{discr_form_FTS} the definitions in \cref{eq:sys2}, and setting $\mathbf{m} = \nu \epsilon \mathbf{\mathfrak{K}}$ and $\mathbf{M}_0 = \mathbf{I} + \mathbf{m} \Delta t$ (with $\mathbf{I}$ the identity matrix), after simple manipulations, \cref{eq:discr_f_1,eq:discr_f_2,eq:discr_f_3} become
\begin{subequations} \label{eq:ppcp}
	\begin{gather} 
		\mathbf{M}_0\frac{\left(\mathbf{u}^{k+1/3}-\mathbf{u}^k\right)}{\Delta t} + \mathbf{u}^{k+1/3} \cdot\ \nabla \left(\frac{\mathbf{u}^{k+1/3}}{\epsilon}\right) + \epsilon \nabla \Psi^k - \nu \nabla^2\mathbf{u}^{k-1/3} +\mathbf{m}\mathbf{u}^k = 0,  \label{eq:Ps1}\\
		\mathbf{M}_0 \frac{\left(\mathbf{u}^{k+2/3}-\mathbf{u}^{k+1/3}\right)}{\Delta t} - \nu \nabla^2\mathbf{u}^{k+2/3} + \nu \nabla^2\mathbf{u}^{k-1/3} = 0,   \label{eq:cs1}\\
		\mathbf{M}_0 \frac{\left(\mathbf{u}^{k+1}-\mathbf{u}^{k+2/3}\right)}{\Delta t}+\epsilon \nabla \Psi^{k+1} -\epsilon \nabla \Psi^k = 0. \label{eq:cs2}
	\end{gather}
\end{subequations}

We discretize the domain $\Omega$ using unstructured triangulations $\Omega_T$ of $N_T$ non-overlapping triangles and $N$ nodes. The computational mesh satisfies the extended Delaunay property as defined in \cite{30} (see Fig. 1 in the referred paper), which can always be obtained in the 2D case (see \cite{31} and literature therein). The reason why we use Delaunay meshes will be explained in the following Sections. A triangle $e$ is called the (computational) cell or (computational) element and the triangle side is called the (element) interface. \par
Inside each triangle $e$ the velocity vector $\mathbf{u}_e$ is assumed $\mathbf{u}_e \in \mathfrak{R}_e$, where $\mathfrak{R}_e$ is the lowest-order Raviart-Thomas ({\bf RT0}) space function \cite{24}, whose basic properties are briefly summarized in Appendix \ref{App1}. Thanks to the {\bf RT0} properties, the velocity components are piecewise constant inside each triangle $e$ if $\sum_{j=1}^3Q_j^e=0$ (where $Q_j^e$ is the normal flux crossing side $j$ of $e$, positive outward, i.e., one of the three DOFs of the $\mathfrak{R}_e$ space). If this condition is satisfied, $\nabla \cdot \mathbf{u}_e=0$  $\forall \mathbf{x}\in e$, $\forall e \in \Omega_T$, and, if the normal fluxes of two neighboring triangles are equal in value and opposite in sign along the common side, both local and global mass continuity are preserved. \par
The kinematic pressure $\Psi$ is assumed to be piecewise linear inside each triangle $e$ according to the nodal values, as explained in \cite{29, 30}. \par
In the present paper, we specifically adapt the procedure proposed in \cite{29, 30} to account for the modified system in \cref{eq:governing_Eqq} of governing equations, compared to the classical Navier-Stokes equations. \par
\begin{figure}[h!]
	\centering
	\includegraphics[width=0.5\textwidth]{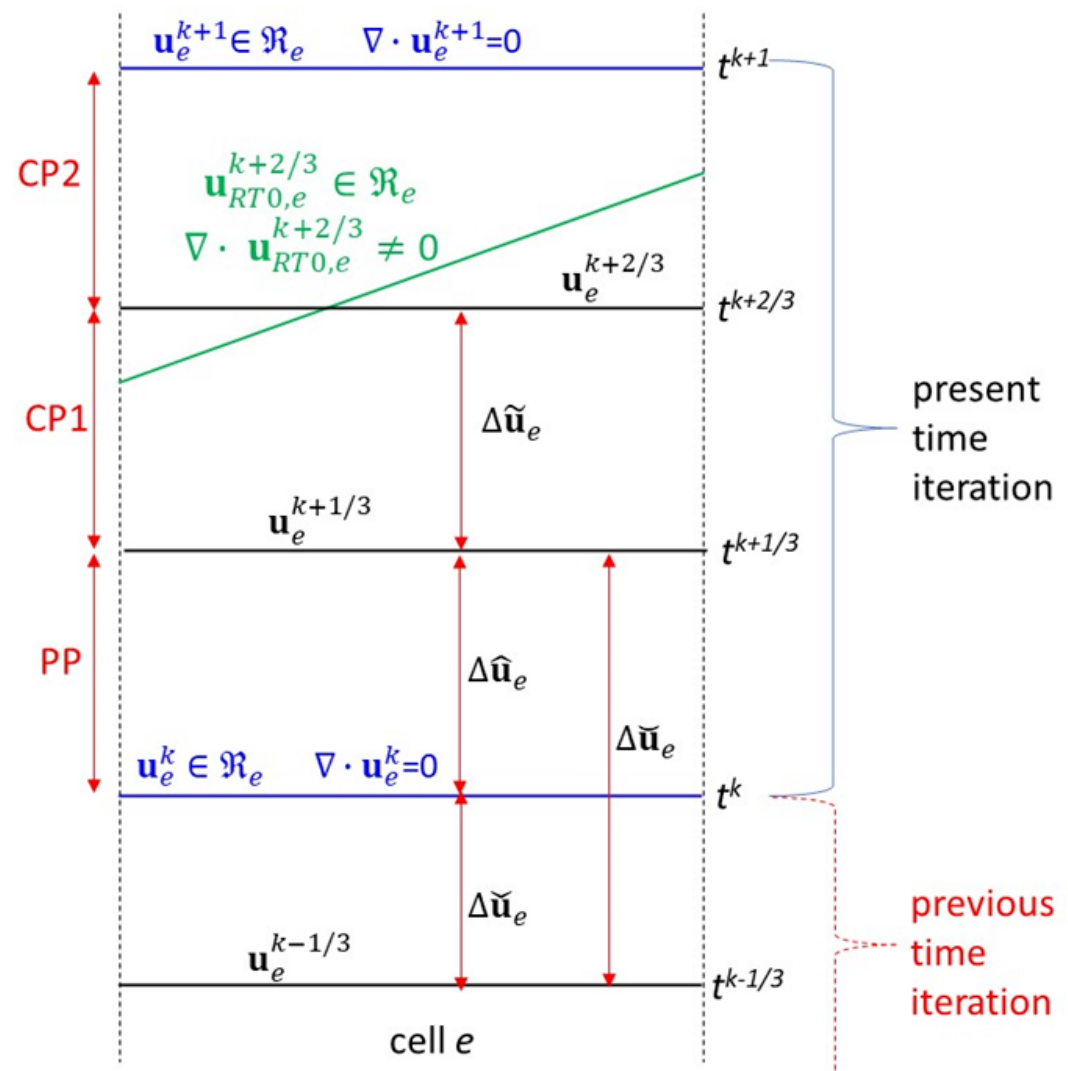}
	\caption{1D sketch of the velocity update within cell $e$ during the algorithm steps for the current time iteration}
	\label{fig_summ_velo}
\end{figure}

	
Below we give the outline of the proposed algorithm, whose numerical details are reported in \cref{algo_details}. In \cref{fig_summ_velo} we show the sequence of the algorithm steps, along with the associated representation of the velocity field within any cell $e \in \Omega_T$. \par	
\begin{itemize}
	\item Beginning of each time iteration (time level $t^k$). $\mathbf{u}_e^k \in \mathfrak{R}_e$, the corresponding normal fluxes are continuous along each element interface and $\nabla \cdot \mathbf{u}_e=0$, $\forall \mathbf{x}\in e$, $\forall e \in \Omega_T$. 
	\item Solution of the PP (from time level $t^k$ to time level  $t^{k+1/3}$). After integration in space over each mesh element of \cref{eq:Ps1}, PP is solved in its time integral form. A ``local element update'' of the velocity field is performed, by computing a piecewise constant correction $\Delta \hat{\mathbf{u}}$ of $\mathbf{u}_e$ within each $e$ (see \cref{fig_summ_velo}). This disrupts the continuity of the normal fluxes (DOFs of $\mathfrak{R}_e$) at each element interface. At the end of PP (time level $t^{k+1/3}$) $\mathbf{u}^{k+1/3}$ is piecewise constant, but does not satisfy local and global mass balance. Numerical details of this prediction problem are presented in \cref{PS}. 
	\item  Solution of the CP1 (from time level $t^{k+1/3}$ to time level  $t^{k+2/3}$). After integration in space over each mesh element, we solve \cref{eq:cs1} in its differential form. Starting from the solution $\mathbf{u}_e^{k+1/3}$, we perform, $\forall e \in \Omega_T$, a ``local element update'' of the velocity field by computing a second piecewise constant correction $\Delta \tilde{\mathbf{u}}$ of $\mathbf{u}_e$ within each $e$ (see \cref{fig_summ_velo}). As for the PP, at the end of CP1 (time level $t^{k+2/3}$) $\mathbf{u}^{k+2/3}$ is piecewise constant, but does not satisfy local and global mass balance. After the solution of CP1, before CP2, (time level $t^{k+1/3}$) we re-establish the normal flux continuity at each element interface, by averaging the fluxes computed according to the velocity field $\mathbf{u}^{k+2/3}$ in the two elements sharing the side. The averaged normal fluxes do not satisfy mass balance since $\mathbf{u}^{k+2/3}$ has been obtained by solving, during the PP and CP1, momentum equations only. The  $\mathbf{RT0}$ velocity vector, $\mathfrak{u}_{RT0}^{k+2/3}$,  associated to the continuous normal fluxes is piecewise-linear within each cell but not divergence free (see \cref{fig_summ_velo}). Numerical details of the CP1 are given in \cref{CS1}.
	\item Solution of the CP2 (from time level $t^{k+2/3}$ to time level  $t^{k+1}$). After integration in space over each mesh element, we solve \cref{eq:cs2} in its differential form. We re-establish local and global mass balance by adding corrective normal fluxes to those computed at the end of CP1. These corrective fluxes are calculated to impose the divergence free condition of the final velocity $\mathbf{u}^{k+1}$. After the solution of CP2, $\mathbf{u}_e^{k+1} \in \mathfrak{R}_e$ and it is piecewise constant within each cell $e$ (see \cref{fig_summ_velo}). Numerical details of the CP2 are presented in \cref{CS2}.
\end{itemize}

The prediction problem is solved by applying a Finite Volume MArching in Space and Time procedure (MAST). As mentioned in the Introduction, this has already been applied in other contexts. One of the main advantages of this procedure is that it performs a sequential solution of small Ordinary Differential Equations (ODEs) systems, one for each computational cell. This allows an ``explicit handling'' of the non-linear convective inertial momentum terms in \cref{eq:momentum}, (i.e., the convective inertial terms within each cell are updated in the time interval $\left[t^k, t^{k+1/3}\right]$ separately from the other cells) avoiding the solution of large systems with non-symmetric matrices as in other numerical schemes, e.g., \cite{32} (\cite{26, 27, 28, 29, 30}). The MAST procedure has shown numerical stability for Courant-Friedrichs-Levy (CFL) numbers greater than 1 (see \cite{26, 27, 28, 29, 30} and references therein). \par 
In both corrector problems, large linear systems of dimension ${N_T}$ are solved, with sparse, symmetric, and, if the Delaunay mesh property holds, also positive definite matrices. This ensures that the system matrices are $\mathcal{M}$-matrices, which avoids nonphysical oscillations in the numerical solution \cite{33}. We apply a mass lumping procedure, similar to the one proposed in \cite{34} for Mixed Hybrid Finite Element, which is well-suited if the Delaunay mesh property holds \cite{29, 30}. The coefficients of the matrices of CP1 and CP2 are constant in time, which makes computations efficient, since their assembly and factorization are performed only once, before the beginning of the time loop. \par

\subsection{Numerical details of the algorithm steps} \label{algo_details}
In the following Sections, $A_e$ marks the area of element $e$, $ep$ is one of the neighboring elements of $e$, and the common side is marked as $j$ and $jp$ in the local numeration of $e$ and $ep$, respectively ($j$, $jp$ = 1, 2, 3). The symbols $\bar{(\cdot)}_e$ and $\bar{(\cdot)}_j$ denote the spatial average of variable $\left(\cdot\right)$ inside each element $e$ and along side $j$, respectively, computed according to the nodal values of $\left(\cdot\right)$. If $\left(\cdot\right)$ is a tensor, the symbols $\bar{(\cdot)}_e$ and $\bar{(\cdot)}_j$ denote the average values of the tensor coefficients, computed according to the coefficients in the three nodes of triangle $e$ and the two nodes of side side $j$, respectively. 
\subsubsection{Predictor problem} \label{PS}	   
Integrating in space \cref{eq:Ps1} over each element $e \in \Omega_T$ and left-multiplying by matrix $\bar{\mathbf{M}}_0^{-1}$, we obtain
\begin{equation}
	\begin{gathered} \label{eq:MAST_PS_5}
		\frac{\left(\mathbf{u}^{k+1/3}-\mathbf{u}^k\right)}{\Delta t} A_e + \bar{\mathbf{M}}_0^{-1} \int_{A_e}\left(\mathbf{u}^{k+1/3} \cdot\ \nabla \left(\frac{\mathbf{u}^{k+1/3}}{\epsilon}\right)\right)  \textrm{d}A + \\
		\bar{\mathbf{M}}_0^{-1} \int_{A_e} \left(\epsilon \nabla \Psi^k - \nu \nabla^2\mathbf{u}^{k-1/3} +\mathbf{m}\mathbf{u}^k\right) \textrm{d}A = 0
	\end{gathered}.
\end{equation}

According to the MAST approach, the scalar momentum equations along the $x$ and $y$ directions are solved inside each triangle $e$ separately from those of the other cells. This is possible if at the beginning of each time step (time level $t^k$) we perform a sorting operation of all the cells in the domain according to the direction of the velocity vector $\mathbf{u}^k$. This is a fast operation as described in \cite{29}. At the end of the sorting operation, a rank $R_e$ is assigned to each cell, such that $R_e> R_{ep}$, with $R_{ep}$ the rank of any neighboring $ep$ triangles, whose common side shared with $e$ is crossed by a flux entering $e$ from $ep$. \par 
The solution of \cref{eq:MAST_PS_5} depends on the pressure gradient and viscous and drag forces obtained from the previous time step, and, thanks to the sorting operation, on the incoming momentum flux from neighboring cells $ep$ with $R_e>R_{ep}$. These are the reasons why  \cref{eq:MAST_PS_5} can be solved within the interval $\left[0 - \Delta t\right]$ (e.g. from $t^k$ to $t^{k+1/3}$), as a system of two Ordinary Differential Equations (ODEs) for the $u$ and $v$ unknowns, $\forall e \in \Omega_T$ \cite{26, 27, 28, 29, 30}. \par 
As specified in \cref{FTS_procedure}, at time level $t^k$ velocity $\mathbf{u}_e^k \in \mathfrak{R}_e$ is divergence free. $\forall e \in \Omega_T$ we define a piecewise constant velocity vector correction $\Delta \hat{\mathbf{u}}_e$, such that
\begin{equation} 
	\mathbf{u}_e(t)=\mathbf{u}_e^k+\Delta \hat{\mathbf{u}}_e(t)  \quad 0 \le t \le \Delta t  \quad \textrm{with} \quad \Delta\hat{\mathbf{u}}_e(0)=0.
\end{equation}

Applying the Green lemma to the second term on the l.h.s. of \cref{eq:MAST_PS_5}, we rewrite the momentum equilibrium equation as
\begin{equation} \label{eq:MAST_ODE_fw}
	\begin{gathered}
		\frac{d \left(\Delta\hat{\mathbf{u}}_e\right)}{d t}A_e + \sum_{j=1}^3\ \phi_j  \bar{\mathbf{M}}_{0,j}^{-1} \mathfrak{M}_j^{e,out}\left(t\right)\ + \sum_{j=1}^3\ \left(1-\phi_j\right) \bar{\mathbf{M}}_{0,j}^{-1} \bar{\mathfrak{M}}_{0,j}^{e,in}\left(t\right)\ + \\
		\bar{\mathbf{M}}_{0,e}^{-1}  \left(\mathbf{\mathfrak{S}}_{\Psi,e}^{k}+\mathbf{\mathfrak{V}}_e^{k-1/3} +\bar{\mathbf{m}}_e \mathbf{u}_{e}^k \right)A_e=0
	\end{gathered},
\end{equation}
\noindent where $\mathfrak{M}_j^{e,out}$ is the leaving momentum flux from side $j$ of $e$ to $ep$ with $R_{ep}>R_e$,
\begin{equation}
	\mathfrak{M}_j^{e,out}\left(t\right)\ = l_j^e \frac{\mathbf{u}_e\left(t\right)\ }{\bar{\epsilon}_j} \max \left(0,\mathbf{u}_e\left(t\right)\ \cdot \mathbf{n}_j\right),
\end{equation}
\noindent and $l_j^e$ is the length of side $j$. $\bar{\mathfrak{M}}_j^{e,in}$ in \cref{eq:MAST_ODE_fw} is the mean in time value of the incoming momentum flux crossing side $j$, oriented from $ep$ to $e$, with $R_{ep}<R_e$, known from the solution of the previously solved cells and computed as specified in Section 3.1 in \cite{29}, $\phi_j=1$ if $e$ shares its side $j$ with triangles $ep$ with $R_{ep}>R_e$, or if it is a  boundary side with leaving momentum flux, in the opposite case $\phi_j=0$. $\mathbf{\mathfrak{V}}_e^{k-1/3}$ and $\mathbf{\mathfrak{S}}_{\Psi,e}^k$ are the sum of viscous and kinematic pressure forces, respectively, computed in cell $e$ during the previous time step, as better specified in \cref{CS1} and \cref{CS2}. \par
\cref{eq:MAST_ODE_fw} is called \textit{MAST forward step}. ODEs systems in \cref{eq:MAST_ODE_fw} are sequentially solved, one for each cell. We apply a Runge-Kutta method with adjustable time-step size within the interval $\left[0, \Delta t\right]$ \cite{35}. We proceed, in the sequential solution, from the triangles with the smallest rank to the triangles with the highest rank. The direction of $\mathbf{u}_e$ could change within $\left[0 - \Delta t\right]$ during the solution of system in \cref{eq:MAST_ODE_fw}, and we could compute momentum fluxes going from $e$ to $ep$ with $R_{ep}<R_e$. These momentum fluxes are neglected during the \textit{MAST forward step}. To restore the momentum balance, after the \textit{MAST forward step}, we perform a \textit{MAST backward step} proceeding from the cells with the highest rank to the cells with the lowest rank. During the \textit{MAST backward step} only the inertial terms are retained,
\begin{equation} \label{eq:MAST_ODE_bw}
	\frac{d \left(\Delta\hat{\mathbf{u}}_e\right)}{d t}A_e +\sum_{j=1}^3\ \left(1-\phi_j \right) \bar{\mathbf{M}}_{0,j}^{-1}  \mathfrak{M}_j^{e,out}\left(t\right)\ + \sum_{j=1}^3\ \phi_j  \bar{\mathbf{M}}_{0,j}^{-1} \bar{\mathfrak{M}}_j^{e,in}\left(t\right) = 0.
\end{equation}

The initial solution of the \textit{MAST backward step} is the final solution computed at the end of the \textit{MAST forward step}. The boundary conditions of the predictor problem are assigned as specified in \cite{29}. As mentioned in \cref{FTS_procedure}, the ``local element update'' by the computation of $\Delta\hat{\mathbf{u}}_e$ $\forall e \in \Omega_T$ disrupts the continuity of the normal fluxes at each element interface, and velocity $\mathbf{u}^{k+1/3}$ is not divergence free.	

\subsubsection{$1^{st}$ corrector problem} \label{CS1}   
Integrating \cref{eq:cs1} in space over each element $e \in \Omega_T$ we get 
\begin{equation} \label{eq:CS1_new}
	\bar{\mathbf{M}}_e   \left(\textbf{u}^{k+2/3}-\textbf{u}^{k+1/3}\right) A_e= \nu \int_{A_e}  \nabla^2 \mathbf{u}^{k+2/3} \textrm{d}A_e  - \nu \int_{A_e} \nabla^2 \mathbf{u}^{k-1/3} \textrm{d}A_e,
\end{equation}
\noindent where matrix $\mathbf{M}=\mathbf{M}_0/\Delta t$, with matrix $\mathbf{M}_0$ defined in \cref{FTS_procedure}. Applying the Green lemma to the integrals on the r.h.s. of \cref{eq:CS1_new}, we get \cref{eq:CS1_2}, which forms a system to solve from $t^{k+1/3}$ to $t^{k+2/3}$ for the $u$ and $v$ unknowns $\forall e \in \Omega_T$,
\begin{equation} \label{eq:CS1_2}
	\bar{\mathbf{M}}_e \left(\textbf{u}^{k+2/3}-\textbf{u}^{k+1/3}\right)A_e =\nu \oint_{L_e} \ \frac{\partial \mathbf{u}_e^{k+2/3}}{\partial n} dl\ - \nu \oint_{L_{e}} \ \frac{\partial \mathbf{u}_{e}^{k-1/3}}{\partial n} dl .
\end{equation}
$\frac{\partial \textbf{u}}{\partial n}$ is the derivative of the velocity vector along the orthogonal direction to the boundary $L_e$ of triangle $e$ and $\oint_{L_e}$ represents the integral over the three sides of triangle $e$. As mentioned at the end of \cref{FTS_procedure}, inside each triangle $e$, we apply a mass-lumping Mixed Hybrid Finite Element procedure to solve system in \cref{eq:CS1_2}, setting (see \cite{29, 30} and literature therein and \cite{34})
\begin{equation} \label{eq:cs1_discr}
	\nu \oint_{L_e} \ \frac{\partial \mathbf{u}_e^{\ast}}{\partial n} dl \ = \nu \sum_{j=1}^3 \ \frac{\left(\textbf{u}_e^{\ast}-\textbf{u}_{ep}^{\ast}\right)}   {d_{e,ep}} l_j^e\ \qquad \ast=\begin{cases}
		k+2/3 \\ k-1/3
	\end{cases},
\end{equation}
\noindent where $d_{e,ep}=\left|\mathbf{x}_e-\mathbf{x}_{ep}\right| sign$, $\left|\mathbf{x}_e-\mathbf{x}_{ep}\right|$ is the distance between the circumcenters of triangles $e$ and $ep$ and $sign$ = 1 or -1 depending on whether or not the mesh satisfies the Delaunay property \cite{29, 30}. The other symbols have been previously specified. $\forall e \in \Omega_T$ we introduce two new piecewise-contant vectors $\Delta \tilde{\mathbf{u}}_e$ and $\Delta \breve{\mathbf{u}}_e$, such that
\begin{equation} \label{eq:cs1_velo}
	\textbf{u}_e^{k+2/3}-\textbf{u}_e^{k+1/3}=\Delta \tilde{\mathbf{u}}_e \qquad \textbf{u}_e^{k+1/3}-\textbf{u}_e^{k-1/3}=\Delta \breve{\mathbf{u}}_e,
\end{equation}
\noindent where $\Delta \breve{\mathbf{u}}_e = \Delta \hat{\mathbf{u}}_e +\Delta \check{\mathbf{u}}_e$, with $\Delta \check{\mathbf{u}}_e = \textbf{u}_e^k-\textbf{u}_e^{k-1/3}$ is known. After some manipulations, system in \cref{eq:CS1_2} can be written as	
\begin{equation} \label{eq:cs1_sys}
	\bar{\mathbf{M}}_e  \Delta \tilde{\mathbf{u}}_eA_e = \nu\sum_{j=1}^3 \ \frac{\Delta \tilde{\mathbf{u}}_e-\Delta \tilde{\mathbf{u}}_{ep}}{d_{e,ep}} l_j^e - \nu \sum_{j=1}^3 \frac{\Delta \breve{\mathbf{u}}_e-\Delta \breve{\mathbf{u}}_{ep}}{d_{e,ep}} l_j^e.
\end{equation}
We solve one system in \cref{eq:cs1_sys} for each component of $\Delta \tilde{\mathbf{u}}_e$ unknown, $\Delta \tilde{u}_e$ and $\Delta \tilde{v}_e$, respectively, 
\begin{subequations} \label{eq:cs1_sys2}
	\begin{gather}
		\begin{cases} \label{eq:cs1_sys2_2}
			\bar{\mathbf{M}}_{x,e} \begin{pmatrix}
				\Delta \tilde{u}_e \\ \Delta \tilde{v}_e
			\end{pmatrix}A_e= \nu \sum_{j=1}^3 \ \frac{\Delta \tilde{u}_e-\Delta \tilde{u}_{ep}}{d_{e,ep}} l_j^e - \nu \sum_{j=1}^3 \frac{\Delta \breve{u}_e-\Delta \breve{u}_{ep}} {d_{e,ep}} l_j^e \\ \bar{\mathbf{M}}_{y,e} \begin{pmatrix}
				\Delta \tilde{u}_e  \\ \Delta \tilde{v}_e
			\end{pmatrix}A_e = \nu \sum_{j=1}^3 \ \frac{\Delta \tilde{v}_e-\Delta \tilde{v}_{ep}}{d_{e,ep}} l_j^e- \nu \sum_{j=1}^3 \frac{\Delta \breve{v}_e-\Delta \breve{v}_{ep}} {d_{e,ep}} l_j^e
		\end{cases}, \\
		\qquad \text{with} \qquad \bar{\mathbf{M}}_{x,e}=\begin{pmatrix} 
			\bar{M}_{1,1}^e\\\bar{M}_{1,2}^e
		\end{pmatrix}^T \qquad \bar{\mathbf{M}}_{y,e}=\begin{pmatrix}
			\bar{M}_{2,1}^e\\\bar{M}_{2,2}^e
		\end{pmatrix}^T,
	\end{gather}	 
\end{subequations}
\noindent and we apply the iterative procedure described in Appendix \ref{App2}. Typically three/four iterations are enough to satisfy the convergence of the iterative procedure and the computational effort required for solving the CP1 step is very small compared to the other algorithm steps, as shown in \cref{test1}. BCs of CP1 are set as in \cite{29}. After solving the systems in \cref{eq:cs1_sys2_2}, we update the velocity at time level $t^{k+2/3}$ according to the first relationship in \cref{eq:cs1_velo}. Due to the ``local element update'' operation performed during the CP1, normal flux continuity at each element interface is not yet recovered and $\mathbf{u}_e^{k+2/3}$ is not divergence free. \par
The sum of the viscous forces $\mathbf{\mathfrak{V}}_e^{k-1/3}$ for the next time iteration in \cref{eq:MAST_ODE_fw} is computed as
\begin{equation}
	\mathbf{\mathfrak{V}}_e^{k-1/3} =  \nu \sum_{j=1}^3\ \frac{\left(\textbf{u}_e^{k+2/3}-\textbf{u}_{ep}^{k+2/3}\right)}{d_{e,ep}} l_j^e .
\end{equation}
At the end of CP1, (time level $t^{k+2/3}$) we compute $\forall e \in \Omega_T$ the velocity vector $\mathbf{u}_{RT0,e}^{k+2/3} \in \mathfrak{R}_e$ according to \cref{eq:RT0}, where the normal flux $Q_j^e$ is given in \cref{eq:flux_k_+_2/3}
\begin{equation} \label{eq:flux_k_+_2/3}
	Q_j^e=\bar{Fl}_j^e=\frac{ \left(\mathbf{u}_e^{k+2/3} \cdot \mathbf{n}_j^e \right)A_{ep} - \left(\mathbf{u}_{ep}^{k+2/3} \cdot \mathbf{n}_{jp}^{ep}\right)A_e}{\left(A_e + A_{ep} \right)} l_j^e \qquad j, jp = 1,2,3 ,
\end{equation} 
\noindent where $\bar{Fl}_j^e$ is the weighted mean flux crossing side $j$ of $e$ computed according to the fluxes $\mathbf{u}_e^{k+2/3} \cdot \mathbf{n}_j^e l_j^e$ and $\mathbf{u}_{ep}^{k+2/3} \cdot \mathbf{n}_{ep}^{k+2/3} l_j^e$ crossing the common side shared by $e$ and $ep$ at time level $t^{k+2/3}$. According to \cref{eq:flux_k_+_2/3}, the fluxes $\bar{Fl}_j^e$ are continuous, i.e, $\bar{Fl}_j^e = -\bar{Fl}_{jp}^{ep}$, but do not satisfy mass balance, i.e., $\sum_{j=1}^3 \bar{Fl}_j^e \neq 0$. This is because so far, the velocity field has been updated from the initial state $\mathbf{u}_e^k$, by solving, in the PP and CP1, momentum equilibrium equations only. This implies that $\nabla \cdot \mathbf{u}_{RT0,e}^{k+2/3} \neq 0$ and velocity $\mathbf{u}_{RT0,e}^{k+2/3}$ is piecewise linear inside $e$. The velocity $\mathbf{u}_{RT0,e}^{k+2/3}$, as well as the continuous fluxes $\bar{Fl}_j^e$ are used for the solution of the $2^{nd}$ correction problem, as explained in the following Section.\par	
The diagonal and off-diagonal matrix coefficients of system in \cref{eq:cs1_sys2_2} $\mathcal{A}_{e,e}^{x\left(y\right),CP1}$ and $\mathcal{A}_{e,ep}^{x\left(y\right),CP1}$, as well as the coefficients $\mathcal{B}_e^{CP1}$ of the source term vector are given in  \cref{eq:CS1_sys_coeff}
\begin{equation} \label{eq:CS1_sys_coeff}
	\begin{gathered}
		\mathcal{A}_{e,e}^{x,CP1}\,=\,\bar{M}_{1,1}^e A_e+ \nu \sum_{j=1}^3 \frac{l_j^e}{d_{e,ep}} \qquad  \mathcal{A}_{e,e}^{y,CP1}=\bar{M}_{2,2}^e A_e + \nu \sum_{j=1}^3 \frac{l_j^e}{d_{e,ep}}, \\ 
		\mathcal{A}_{e,ep}^{x(y),CP1}=-\nu \frac{l_j^e}{d_{e,ep}}, \\
		\mathcal{B}_e^{CP1}=\nu \sum_{j=1}^3 \frac{\Delta \breve{u}_e-\Delta \breve{u}_{ep}} {d_{e,ep}} l_j^e -\bar{M}_{1,2}^e A_e \Delta \tilde{v}_e^i, \\
		\mathcal{B}_e^{CP1}=\nu \sum_{j=1}^3 \frac{\Delta \breve{v}_e-\Delta \breve{v}_{ep}} {d_{e,ep}} l_j^e -\bar{M}_{2,1}^e A_e \Delta \tilde{u}_e^i.		
	\end{gathered}	
\end{equation}
The matrices of systems in \cref{eq:cs1_sys2_2} are sparse and symmetric. If the mesh satisfies the Dalaunay property, $d_{e,ep}$ are positive, and the matrices are also positive-definite, so that the $\mathcal{M}$-matrix property is guaranteed \cite{33}. The systems in \cref{eq:cs1_sys2_2} are solved by a fast and efficient Preconditioned Conjugate Gradient (PCG) method \cite{36,37} with an incomplete Cholesky factorization \cite{38,39}. The matrix coefficients only depend on $\nu$, $M_{i,j}^e$ and geometrical quantities, so that the matrices of systems in \cref{eq:cs1_sys2_2} are only factorized once, before the beginning of the time loop, saving a lot of computational effort. \par 

\subsubsection{$2^{nd}$ corrector problem} \label{CS2}   
\cref{eq:cs2} is rewritten as 
\begin{equation} \label{eq:Cs2_1}
	\mathbf{M}  \left(\mathbf{u}^{k+1}-\mathbf{u}^{k+2/3}\right) +\epsilon \left(\nabla \Psi^{k+1}-\nabla \Psi^k\right) = 0,
\end{equation}
\noindent with matrix $\mathbf{M}$ defined in \cref{CS1}. Introducing the scalar variable $\eta$ such that
\begin{equation} \label{eq:Cs2_2}
	\mathbf{M} \left(\mathbf{u}^{k+1}-\mathbf{u}_{RT0}^{k+2/3}\right)=\epsilon  \nabla \eta,
\end{equation}
\noindent and left-multiplying both sides of \cref{eq:Cs2_2} by $\mathbf{M}^{-1}$, we obtain
\begin{equation} \label{eq:Cs2_3}
	\left(\mathbf{u}^{k+1}-\mathbf{u}_{RT0}^{k+2/3}\right)=\epsilon \mathbf{M}^{-1} \nabla \eta ,
\end{equation}
\noindent with $\mathbf{u}_{RT0}^{k+2/3}$ defined as in \cref{CS1}. Taking the divergence of \cref{eq:Cs2_3} we get
\begin{equation} \label{eq:Cs2_4}
	\nabla \cdot \left(\mathbf{u}^{k+1}-\mathbf{u}_{RT0}^{k+2/3}\right)= \nabla \cdot \left( \mathbf{\Xi} \nabla \eta\right) \quad \textrm{with} \quad \mathbf{\Xi}=\epsilon \mathbf{M}^{-1}.
\end{equation}

\begin{figure}[h!]
	\centering
	\includegraphics[width=0.5\textwidth]{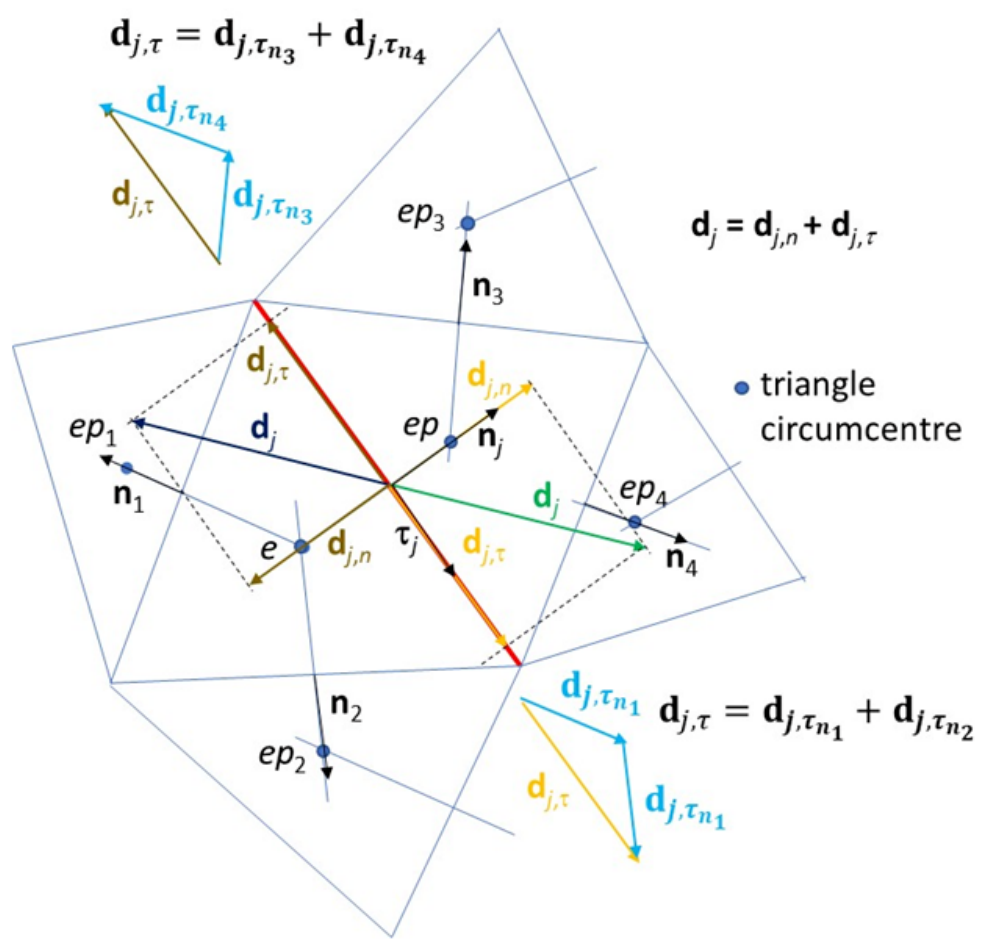}
	\caption{CP2, calculation of $\left( \mathbf{\Xi}_j \mathbf{n}_j \right) \cdot \nabla \eta$}
	\label{fig_cs2}
\end{figure}

Setting $\nabla \cdot \mathbf{u}^{k+1} = 0$, integrating in space over each element $e$ and applying the Green lemma, we obtain from \cref{eq:Cs2_4}
\begin{equation} \label{eq:Cs2_5}
	\sum_{j=1}^3 \left(\bar{Fl}_j^e \right) + \oint_{L_{e}} \left(\bar{\mathbf{\Xi}} \nabla \eta\right) \cdot \mathbf{n} \: \mathrm{d}l= 0 \quad  \textrm{ with } \oint_{L_{e}} \left(\bar{\mathbf{\Xi}} \nabla \eta\right) \cdot \mathbf{n} \: \mathrm{d}l = \sum_{j=1}^3\left(\bar{\mathbf{\Xi}}_j \nabla \eta\right) \cdot \mathbf{n}_j \: l_j^e ,
\end{equation}
\noindent where the average normal flux $\bar{Fl}_j^e$ crossing side $j$ of cell $e$ has been defined in \cref{eq:flux_k_+_2/3}. Since matrix $\bar{\mathbf{\Xi}}_j$ is symmetric and positive definite, we have $\left(\bar{\mathbf{\Xi}}_j \nabla \eta\right) \cdot \mathbf{n}_j = \left( \bar{\mathbf{\Xi}}_j \mathbf{n}_j \right) \cdot \nabla \eta$, and we set $\bar{\mathbf{\Xi}}_j \mathbf{n}_j = \mathbf{d}_j$. Decomposing vector $\mathbf{d}_j$ along the normal and tangential directions to side $j$ ($\mathbf{n}_j$ and $\boldsymbol{\tau}_j$ in \cref{fig_cs2}), $\mathbf{d}_j=\mathbf{d}_{j,n} + \mathbf{d}_{j,\tau}$, and applying a  co-normal decomposition of vector $\mathbf{d}_{j,\tau}$ along the normal directions to the other two sides of cell $e$ ($\mathbf{n}_1$ and $\mathbf{n}_2$ in \cref{fig_cs2}), after some manipulations whose details are given in Appendix \ref{App3}, \cref{eq:Cs2_5} is discretized as 
\begin{equation}
	\begin{gathered} \label{eq:Cs2_6}
		\sum_{j=1}^3 \left(\bar{Fl}_j^e \right) = \sum_{j=1}^3 \left(\frac{\eta_e - \eta_{ep}}{d_{e,ep}} d_{j,n} \right)l_j^e + \\
		\frac{1}{2} \sum_{j=1}^3 \left(\sum_{l=1,2} \frac{\eta_e - \eta_{ep_l}}{d_{e,ep_l}} d_{\tau_{n_l}} \alpha_l + 
		\sum_{m=3,4} \frac{\eta_{ep_m} - \eta_{ep}}{d_{ep,ep_m}} d_{\tau_{n_m}} \alpha_m \right) l_j^e 
	\end{gathered},
\end{equation}

\noindent which forms a system to be solved for the $\eta$ unknowns. With the help of Appendix \ref{App3} and \cref{fig_cs2}, $d_{j,n}=\mathbf{d}_j \cdot \mathbf{n}_j$, $d_{\tau_{n_{l(m)}}}=\mathbf{d}_{j,\tau} \cdot \mathbf{n}_{l(m)} $ (\textit{l=}1,2, \textit{m}=3,4), $\eta_e$, $\eta_{ep}$, $\eta_{ep_{l(m)}}$ are the $\eta$ unknowns in the circumcenters of cells $e$, $ep$, and $ep_{l(m)}$, respectively, and  $d_{e,ep}$, $d_{e,ep_l}$ and $d_{e,ep_m}$ are the distances of the circumcenters of cells $e$ and $ep$, $e$ and $ep_l$, $ep$ and $ep_m$ respectively, times +1 or -1, depending on whether the mesh satisfies or not the Delaunay property (see also \cref{CS1}). $\alpha_{l(m)}=1$ if $\mathbf{d}_{j,\tau} \cdot \mathbf{n}_{l(m)} > 0$ otherwise $\alpha_{l(m)}=-1$.\par 
Instead of solving system in \cref{eq:Cs2_6}, we proceed as follows. The solution of system in \cref{eq:Cs2_7_1} gives an approximate solution of $\eta$, denoted as $\tilde{\eta}$. With this, the final solution can be obtained by solving system in \cref{eq:Cs2_7_2} 
\begin{subequations} \label{eq:Cs2_7}
	\begin{gather}
		\sum_{j=1}^3 \left(\frac{\tilde{\eta}_e - \tilde{\eta}_{ep}}{d_{e,ep}} d_{j,n}\right) l_j^e  = \sum_{j=1}^3 \left(\bar{Fl}_j^e \right), \label{eq:Cs2_7_1} \\			
		\sum_{j=1}^3 \left(Fl_j^{\eta,e}\right) =  \sum_{j=1}^3 \left(\bar{Fl}_j^e \right), \label{eq:Cs2_7_2}
	\end{gather}
\end{subequations}
\noindent where $Fl_j^{\eta,e}$ is the flux crossing side $j$ of $e$ due to $\mathbf{\Xi} \nabla \eta $,
\begin{equation} \label{eta_flux}
	\begin{gathered}
		Fl_j^{\eta,e} =\int_{l_{j,e}} \left(\mathbf{\Xi} \nabla \eta\right) \cdot \mathbf{n}_j dl_j =  \frac{\eta_e-\eta_{ep}}{d_{e,ep}} d_n l_j^e  + \\
		\frac{1}{2} \left(\sum_{l=1,2} \frac{\tilde{\eta}_e - \tilde{\eta}_{ep_l}}  {d_{e,ep_l}}d_{\tau_{n_l}}\alpha_l + \sum_{m=3,4} \frac{\tilde{\eta}_{ep} - \tilde{\eta}_{ep_m}} {d_{ep,ep_m}}d_{\tau_{n_m}}\alpha_m \right) l_j^e
	\end{gathered} .
\end{equation}

The same spatial discretization as in \cref{eq:cs1_discr}  has been applied. \cref{eta_flux} implies that fluxes $Fl_j^{\eta, e}$, as well as $\bar{FL}_j^e$, are continuous for the two neighbor cells $e$ and $ep$, and \cref{eq:Cs2_5} implies that $\sum_{j=1}^3 \left(Fl_j^{\eta, e} + \bar{FL}_j^e\right) = 0$, $ \forall e \in \Omega_T$. Mass conservation along the three sides and inside each cell $e$, $\forall e \in \Omega_T$ is finally recovered at time level $t^{k+1}$. \par 
After the solution of systems in \cref{eq:Cs2_7}, we obtain the velocity vector at the end of the time step $\forall e \in \Omega_T$, $\mathbf{u}_e^{k+1}$, as in \cref{eq:Cs2_3}, where $\mathbf{\Xi }\nabla\eta$, as well as $\mathbf{u}_{RT0}^{k+2/3}$, is a $\mathbf{RT0}$ function, 
\begin{equation}
	\mathbf{u}_{RT0}^{k+2/3}=\sum_{j=1}^3\left(\mathbf{w}_j^e\bar{FL}_j^e\right) \qquad  \mathbf{\Xi }\nabla \eta=\sum_{j=1}^3 \left(\mathbf{w}_j^e Fl_j^{\eta}\right),
\end{equation}
\noindent and we get 
\begin{equation}
	\mathbf{u}_e^{k+1}=\sum_{j=1}^3 \mathbf{w}_j^e \left(\bar{Fl}_j^e + Fl_j^{\eta}\right).
\end{equation} 

According to the properties of the $\mathbf{RT0}$ functions (see \cref{eq:RT0,eq:RT0_prop}), $\nabla \cdot \mathbf{u}_e^{k+1}=0$ $\forall \mathbf{x}\in e$, $\forall e \in \Omega_T$, and the method is \textit{strongly} conservative (see \cref{intro}). 
Once $\mathbf{u}_e^{k+1}$ is known, we get from \cref{eq:Cs2_1}, term $\epsilon_e \nabla \Psi^{k+1}$, $\forall e \in \Omega_T$, as,
\begin{equation} \label{eq:grad_PSI}
	\bar{\epsilon}_e \nabla \Psi_e^{k+1} = \bar{\epsilon}_e \nabla \Psi_e^k + \bar{\mathbf{M}}_e \left(\mathbf{u}_e^{k+1}-\mathbf{u}_e^{k+2/3}\right),
\end{equation}
\noindent where the symbols have been specified before. We set in \cref{eq:MAST_ODE_fw} of the PP of the next time iteration
\begin{equation}
	\mathfrak{G}_{\eta,e}^k=\bar{\epsilon}_e \nabla \Psi^{k+1}\:A_e.
\end{equation}

The boundary conditions of the CP2 are assigned as specified in \cite{29}. The diagonal and off-diagonal matrix coefficients $\mathcal{A}_{e,e}^{CP2}$ and $\mathcal{A}_{e,ep}^{CP2}$ of the systems in \cref{eq:Cs2_7_1,eq:Cs2_7_2} are
\begin{equation} \label{eq:CS2_sys_coeff}
	\mathcal{A}_{e,e}^{CP2}=\sum_{j=1}^3 \frac{l_j^e}{d_{e,ep}}d_{j,n} \qquad \mathcal{A}_{e,ep}^{CP2}=-\frac{l_j^e}{d_{e,ep}}d_{j,n},
\end{equation}
\noindent and the coefficients of the source term vectors of systems \cref{eq:Cs2_7_1,eq:Cs2_7_2}, $\mathcal{B}_{1,e}^{CSP}$ and $\mathcal{B}_{2,e}^{CP2}$ respectively, are
\begin{equation}
	\begin{gathered}
		\mathcal{B}_{1,e}^{CP2} = \sum_{j=1}^3 \left(\bar{Fl}_j^e \right), \\
		\mathcal{B}_{2,e}^{CP2} = \sum_{j=1}^3 \left(\bar{Fl}_j^e \right) - \\ 
		\frac{1}{2}\sum_{j=1}^3 \left(\sum_{l=1,2} \frac{\tilde{\eta}_e - \tilde{\eta}_{ep_l}}  {d_{e,ep_l}}d_{\tau_{n_l}}\alpha_l + \sum_{m=3,4} \frac{\tilde{\eta}_{ep} - \tilde{\eta}_{ep_m}} {d_{ep,ep_m}}d_{\tau_{n_m}}\alpha_m \right)  l_j^e.
	\end{gathered}
\end{equation}

The matrix $\mathbf{\mathcal{A}}_{CP2}$ is sparse and symmetric, yelding the same beneficial matrix properties as for the CP1 problem if the mesh satisfies the Delaunay property. Therefore, systems in \cref{eq:Cs2_7_1,eq:Cs2_7_2} are solved by applying the same procedure as in \cref{CS1}. Moreover, this again allows to perform factorization of matrix $\mathbf{\mathcal{A}}^{CP2}$ only once.\par
The matrix associated with the original system in \cref{eq:Cs2_6} is not symmetric, and the advantage of splitting system in \cref{eq:Cs2_6} into systems in \cref{eq:Cs2_7_1,eq:Cs2_7_2} is twofold: 1) numerical stability is achieved, since matrix $\mathbf{\mathcal{A}}^{CP2}$ is a $\mathcal{M}$-matrix and 2) there is a fast and efficient solution of the PCG method, compared to the standard GMRES, BiCG, CGSquared, or BiCGStab methods, usually applied for the solution of non-symmetric matrix systems. \par 
Observe in \cref{eta_flux} that flux $Fl_j^{\eta, e}$ depends on the six values of $\eta$ in cells $e$, $ep$ and $ep_{l(m)}$ (\textit{l}=1,2, \textit{m}=3,4). Since the $\eta$ values in cells $ep_{l(m)}$ are assumed to be known in the solution of system in \cref{eq:Cs2_7_2}, $Fl_j^{\eta,e}$ depends only on the two unknowns $\eta_e$ and $\eta_{ep}$. For this reason, due to the splitting strategy operated in \cref{eq:Cs2_7_1,eq:Cs2_7_2}, the flux discretization scheme in \cref{eta_flux} can be regarded, \textit{de facto}, as a Two-Point-Flux-Approximation scheme (TPFA). \par 
According to \cite{29}, within each triangle, $\Psi^{k+1}$ (or $\Psi^k$) is piecewise-linear, while $\eta$ is piecewise-quadratic. Since the kinematic pressure is not needed to update the solution at each time step, we compute the nodal values of $\Psi$ at target simulation times only, as explained in \cite{29, 30}.\par

\section{Numerical tests and analysis of the results} \label{num_tests} 

We present five numerical tests. In the first test we analyze the convergence order of the proposed algorithm and the required computational (CPU) times. In the second test, we compare the solution of the presented algorithm with an analytical one provided in the literature for Stokes flow regime and different geometrical domain configurations. The third, fourth and fifth test are related to ``real-world'' applications. In the third and fourth test we compare our numerical solution with averaged pore-scale results for the Stokes and Navier-Stokes flow regimes. In the last test, we provide a \textit{showcase} with different anisotropy tensors, for different Reynolds numbers and a comparison with the solution provided by a TDA code developed in the framework of the open-source software package $\textrm{DuMu}^{\textrm{x}}$ \cite{40} which applies a Multi-Point-Flux-Approximation (MPFA) scheme for calculating the solution in the porous region. \par 
As ICs for all the presented test cases, we adopt zero velocity and zero kinematic pressure in the domain. \par
A very fast off-line \textit{in-house} procedure is adopted to generate the computational mesh and the input data, and assign the ICs and BCs \cite{30}. The output model results are processed with Paraview \cite{41}. In the presented tests, we neglect gravitational effects.\par  
In the following Sections, $\delta_0$ and $\delta_{TL}$ mark the mesh sizes adopted to discretize the bulk free fluid and porous regions, as well as the the transition layer, respectively, while $d_{TL}$ is the width of the transition layer. 

\subsection{Test 1. Convergence test} \label{test1}  

We assume a 2D domain $\Omega=\left[0,2\right]^2 \textrm{ m}^2$ with an internal square porous region $\Omega_{pm}$  $\left[0.5,1.5\right]^2 \textrm{ m}^2$. In the bulk porous region, the porosity is $\epsilon=\epsilon_0=0.4$, and the anisotropy tensor \textbf{K} is given by \cref{aniso}
\begin{equation}\label{aniso}
	\mathbf{K}=\mathbf{R}\mathbf{A}\mathbf{R}^{-1} \quad \mathbf{A}=\left(\begin{array} {cc} 
		\frac{k}{\beta} & 0 \\
		0 & k \end{array}\right) \quad \mathbf{R}=\left( \begin{array}{cc}
		\textrm{cos}\alpha & -\textrm{sin} \alpha\\
		\textrm{sin}\alpha & \textrm{cos}\alpha
	\end{array}\right),
\end{equation}
\noindent with $\beta = 100$, $\alpha =- \pi/4$, such that $k = \num{1.e-6} \textrm{ m}^2$. Matrix \textbf{K} is symmetric and positive definite. Kinematic fluid viscosity is $\nu=\num{1.5e-5} \textrm{ m}^2/\textrm{s}$. The pressure field is given by	
\begin{equation}	\label{eq:1}
	\Psi = a_1 + b_1 x + c_1 x^2 + a_2 + b_2 y + c_2 y^2,
\end{equation}	
\noindent where the values of the coefficients $a_i$, $b_i$ and $c_i$ $\left(i=1,2\right)$ are listed in \cref{tab1}. The analytical velocity solution is constructed to be divergence free and continuous at the interface between free fluid and porous regions, with velocity components 
\begin{equation}
	\begin{gathered} \label{eq:2}
		u=-(K_{1,1} \dfrac{\partial\Psi}{\partial x} + K_{1,2} \frac{\partial\Psi}{\partial y} ), \\
		v=-(K_{2,1} \dfrac{\partial\Psi}{\partial x} + K_{2,2} \frac{\partial\Psi}{\partial y} ).
	\end{gathered}	
\end{equation}

The Reynolds number is computed as $Re=\frac{\|\mathfrak{u}\|_{max}H_{pm}}{\nu}$, where $\|\mathfrak{u}\|_{max}$ is the maximum value of the velocity vector magnitude and $H_{pm}$ is the side of the porous domain, such that $Re\simeq267$. \par	
We assume $\Gamma_n=\Gamma$, except in the upper right corner, where we assign the Dirichlet condition for $\Psi$ equal to the value given by \cref{eq:1}. \par
\begin{figure}[h!]
	\centering
	\begin{subfigure}{.45\textwidth}
		\centering
		\includegraphics[width=\textwidth]{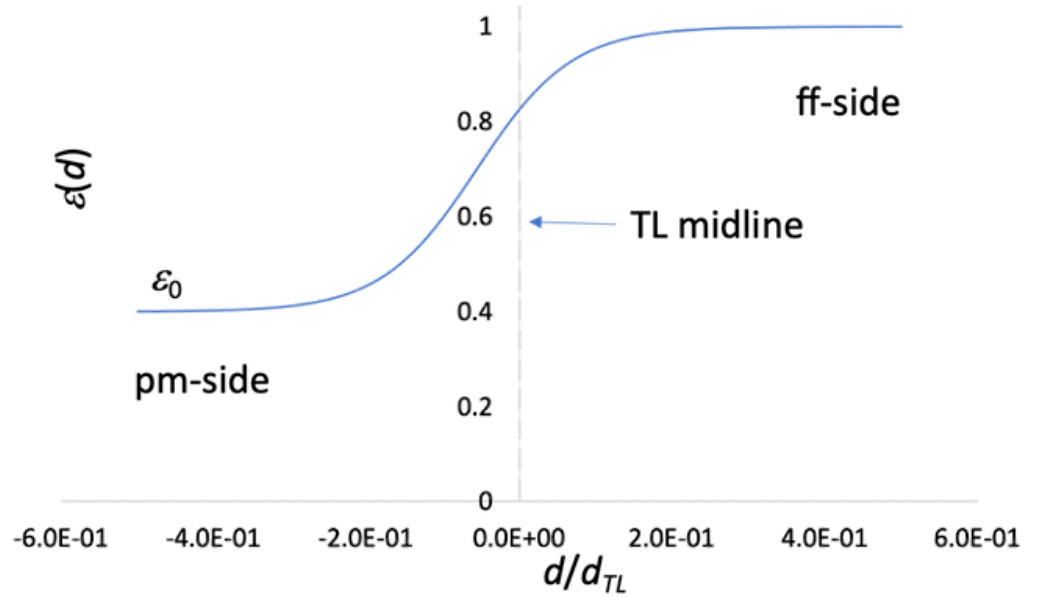}
		\caption{Porosity profile for $d_{\epsilon}<0$}
		\label{fig_poro_profile}%
	\end{subfigure}
	\begin{subfigure}{.45\textwidth} 
		\includegraphics[width=\textwidth]{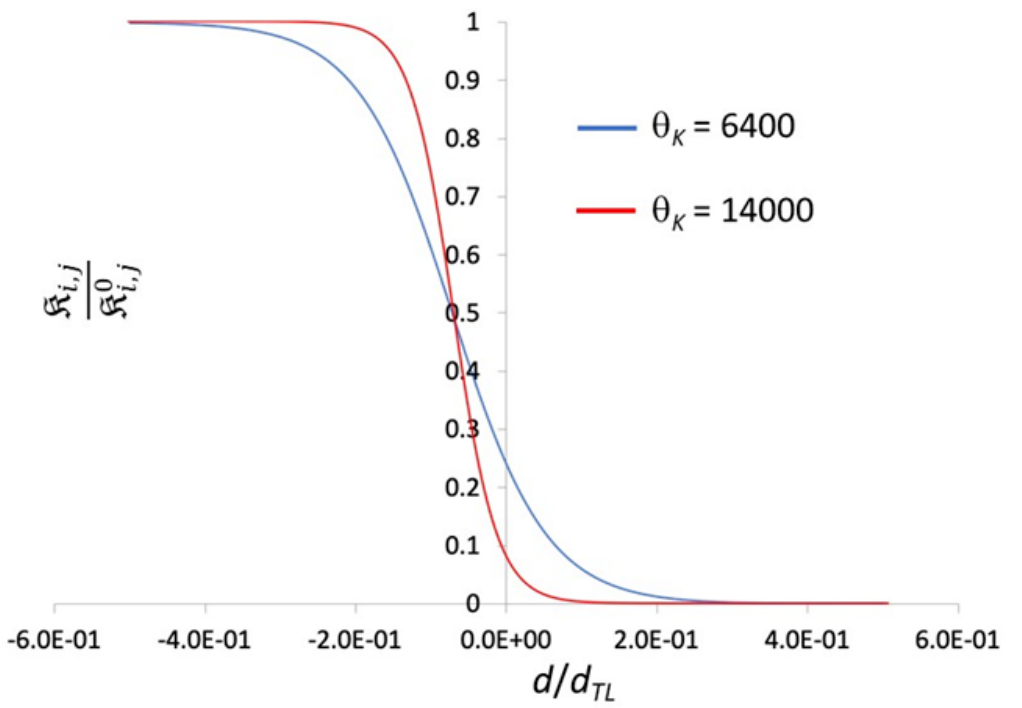}
		\caption{Dimensionless (i,j)-permeability coefficient profiles for $d_K<0$ and two different values of $\theta_K$}
		\label{fig_perm_profiles}%
	\end{subfigure}
	\caption{Porosity and permeability profiles}
	\label{fig_poro_perm_profiles}%
\end{figure}

 
In our scheme, the outer boundary of the transition zone overlaps the outer contour of $\Omega_{pm}$. A continuous variation of the porosity $\epsilon$ and the coefficients $\mathfrak{K}_{i,j}$ of the inverse permeability tensor is assumed within the transition layer as in \cref{eq:3}
\begin{subequations} \label{eq:3}
	\begin{gather} 
		\epsilon(\textbf{x}) = \frac{1}{2}\left(1-\epsilon_0\right)\tanh((d+d_{\epsilon}) \theta_{\epsilon})+\frac{1}{2}\left(1+\epsilon_0\right), \label{eq:3.a}\\
		\mathfrak{K}_{i,j}\left(\textbf{x}\right)= \mathfrak{K}_{i,j}^0 \frac{1}{2}\left(1-\tanh((d+d_k)\theta_K)\right),
		\label{eq:3.b}
	\end{gather}	
\end{subequations}
\noindent where $\mathfrak{K}_{i,j}^0$ are the coefficients of the inverse of the permeability tensor within the bulk $\Omega_{pm}$, computed according to \cref{aniso} and $d=\left|\mathbf{x}_P - \mathbf{x}_{TL,mid}\right| sign$, where $\left|\mathbf{x}_P - \mathbf{x}_{TL,mid}\right|$ is the distance of any point \textit{P} within TL from the center-line of TL and $sign$ = 1 or -1 depending on whether point \textit{P} is located in the half-region of TL on the side of the porous medium or the free fluid region. $d_{\epsilon}$ and $d_K$ are scalar values which control the symmetry of the profiles of $\epsilon$ and $\frac{1}{K_{i,j}}$ with respect to the center-line of TL. Depending on whether they assume positive or negative values, the profiles are shifted towards the free fluid region or to the porous medium region. $\theta_{\epsilon}$ and $\theta_K$ are positive scalar values which control the slope of the profiles (the larger they are, the steeper the profiles). \cref{fig_poro_profile} shows a porosity profile computed setting a negative value of $d_{\epsilon}$, while in \cref{fig_perm_profiles} we plot two dimensionless permeability profiles corresponding to a negative $d_K$ value. \par  
We initially assumed the width of TL to be $d_{TL}=H_{pm}/25$, and the domain is discretize with a coarse mesh ($N_T$ = 16260 triangles and $N$ = 8292), whose maximum values of mesh sizes are  $\delta_0=0.025$ m and $\delta_{TL}=0.02$ m, respectively. Starting from this coarse mesh, we progressively performed four refinement operations by halving $\delta_0=$, leaving $d_{TL}$ unchanged and setting $\delta_0=\delta_{TL}$. The adopted time step size is $\Delta t = \num{2.25e-1}$ s. The maximum $CFL$ value, $CFL_{max}=\|\mathfrak{u}\|_{max} \Delta t/\sqrt{A_e}$ was computed in TL, and $CFL_{max} \simeq 1.387$. At each mesh refinement, we halved $\Delta t$, to avoid increases of $CFL_{max}$. We set both $d_{\epsilon}$ and $d_K$ to zero, and $\theta_{\epsilon}=\theta_K=200$ in \cref{eq:3}.\par	
The $L_2$ norm of the errors of the computed solutions for $u$, $v$ and $\Psi$ with respect to the exact solutions given in \cref{eq:1,eq:2} is computed as  
\begin{equation} \label{error}
	\begin{gathered}
		L_2\left(err\left(q\right)\right) = \sqrt{\sum_{e=1}^{N_T} A_e\left( q_e^n - q_e^{ex}\right)^2 }, \\
		L_2\left(err\left(\Psi\right)\right) = \sqrt{\sum_{i=1}^{N} A_i\left( \Psi_i^n - \Psi_i^{ex}\right)^2},
	\end{gathered}
\end{equation}
\noindent where $q=u \textrm{ or }v$, $A_i$ is the area of the Voronoi polygon associated with node $i$ and superscripts $n$ and $ex$ mark the numerical and exact solutions, respectively. $q_e^n$ and $q_e^{ex}$ are computed in the circumcentre of each cell, while $\Psi_i^n$ and $\Psi_i^{ex}$ are computed in the mesh nodes. If we call $h_l$ the mesh size associated with the $l^{th}$ refinement, and assume that the error associated with the $l^{th}$ refinement is  proportional to a power $r_c$ of $h_l$, we compute the spatial rate of convergence $r_c$ by comparing the errors obtained for the two meshes with the two consecutive linear sizes $h_l$ and $h_{l+1}$,
\begin{equation} \label{conv_order}
	r_c=\ln\left[\frac{L_2\left(err_l\right)}{L_2\left(err_{l+1}\right)}\right] / \ln\left[\frac{h_l}{h_{l+1}}\right].
\end{equation}

In \cref{tab2} we list the $L_2$ norms of the errors and the convergence order $r_c$ of the velocity components is close to 1, due to the piecewise constant approximation of the velocity inside each triangle. $r_c$ of $\Psi$ is smaller than 2. The reason could be that, due to the lack of a specific equation for $\Psi$ in the governing equations, this is indirectly computed from the pressure gradients, as described at the end of  \cref{CS2}.\par 
We also investigate how the size of TL, $d_{TL}$, affects the computed results compared to the exact ones. Starting from the $2^{nd}$ mesh refinement level, we progressively halved $d_{TL}$ as well as its mesh size $\delta_{TL}$, without changing $\delta_0$ in the free fluid and the bulk porous regions. At each refinement of $d_{TL}$, we also halved $\Delta t$ for the aforementioned reasons. Since the assigned analytical solution of the velocity vector depends on the values of the permeability coefficients in the bulk porous region, without any transition of their values close to the interface with the fluid region (see \cref{eq:3}), we expect the numerical solution to get closer and closer to the exact one by refining $d_{TL}$. This is confirmed by the results in \cref{tab_test1_2}. \par
\begin{figure}[h!]
	\centering
	\includegraphics[width=0.7\textwidth]{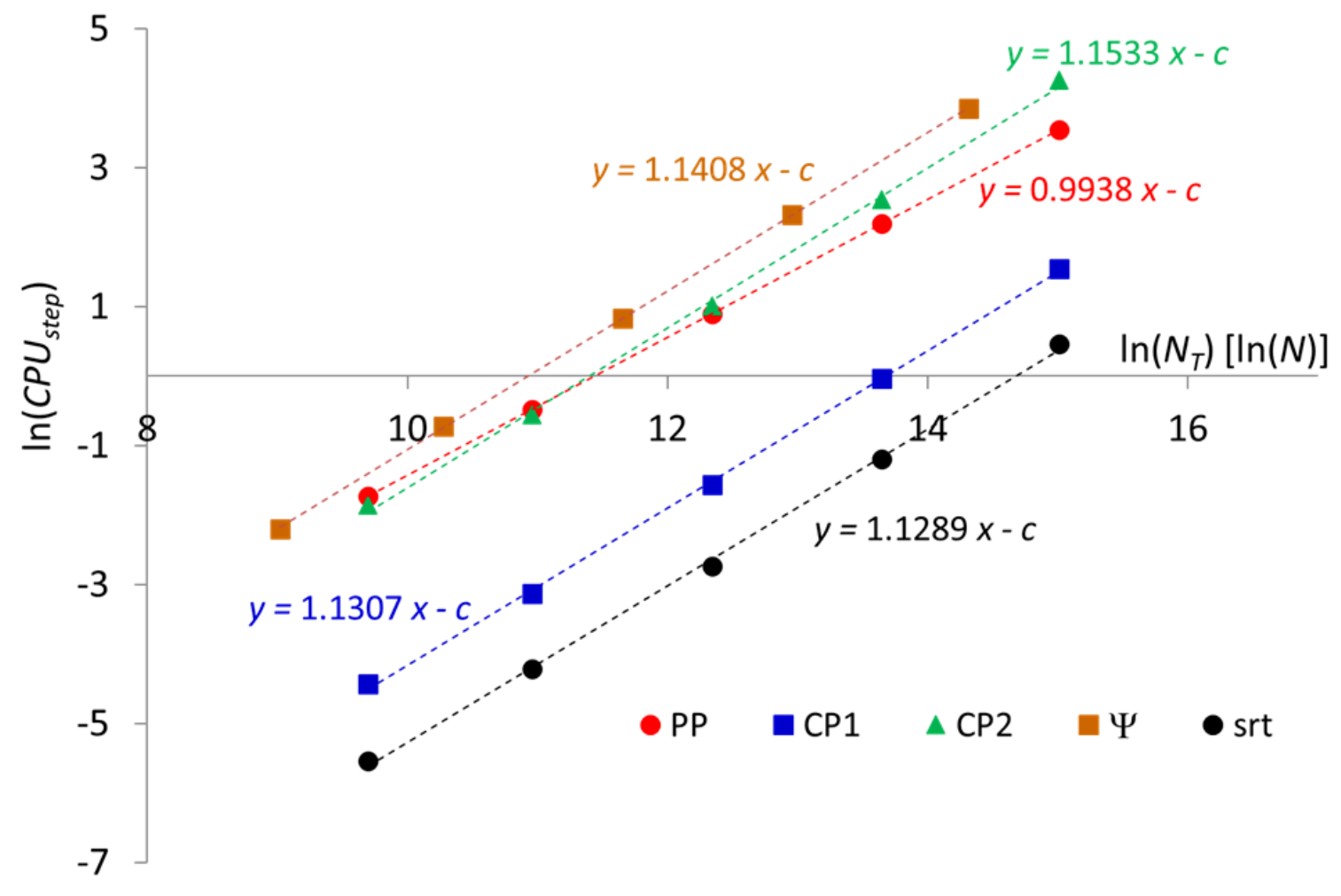}
	\caption{Study of the CPU times of the algorithm steps}
	\label{fig_cpu_times}
\end{figure}

 
We also investigated the computational (CPU) times required by the algorithm steps, PP, CP1, CP2 and sorting cell operation, $CPU_{PP}$, $CPU_{CP1}$, $CPU_{CP2}$, $CPU_{\Psi}$ and $CPU_{srt}$, respectively. Given two real scalar numbers  $c$ and $\omega$, we express the mean value of $\overline{CPU}_{step}$ per time iteration as   	
\begin{equation} \label{CPU_times}
	\overline{CPU}_{step} = \exp\left(c\right) N_T^\omega \qquad \textrm{or} \qquad \overline{CPU}_{step} = \exp\left(c\right) N^\omega.
\end{equation} 
\noindent where ``\textit{step}'' in \cref{CPU_times} corresponds to \textit{PP}, or \textit{CP1}, or \textit{CP2}, or $\Psi$, or \textit{srt}. A single Intel(R) Core(TM) i7-9700K processor at 3.40 GHz was used for the simulation runs. In \cref{fig_cpu_times} we show the computational times in bi-logarithmic scales. Due to the \textit{explicit nature} of the predictor step (i.e., the sequential solution of the ODEs systems during the MAST forward and backward steps), its growth with  $N_T$ is almost linear (the $\gamma$ exponent is slightly smaller than 1). Since the two corrector steps, as well as the computation of the kinematic pressure require the solution of linear systems, their growth with $N_T$ (or $N$) is more than linearly proportional (the associated exponents $\gamma$ are slightly higher than 1, ranging from 1.1307 to 1.1533). $\overline{CPU}_{CP1}$ and $\overline{CPU}_{srt}$ are approximately 1 and 2-3 magnitude orders smaller than $\overline{CPU}_{PP}$ and $\overline{CPU}_{CP2}$. 

\begin{table} [ht] \centering
	\caption{Test 1. Parameter values in \cref{eq:1}}  \label{tab1}
	\begin{tabular}{c|c|c|c|c|c}
		$a_1$ & $a_2$ & $b_1$ & $b_2$ & $c_1$ & $c_2$\\
		\hline
		\num{5e-4} & \num{1.5e-3} & \num{-3.2e-3} & \num{-8e-4} & \num{-1.95e-3} & \num{4e-3}\\
	\end{tabular}
\end{table}

\begin{table}[ht] \footnotesize
	\caption{Test1. $L_2$ norms of errors and $r_c$ for mesh refinement} \label{tab2} 
	\centering 
	\begin{tabular} {c|c|c|c|c|c|c|c}
		\multicolumn{2}{c|}{mesh details} & \multicolumn{2}{c|}{$u$} & \multicolumn{2}{c}{$v$} & \multicolumn{2}{|c}{$ \Psi $} \\
		\hline
		$N$ & $N_T$ & $ L_2 $ & $ r_c $ & $ L_2 $ & $ r_c $ & $ L_2 $ & $ r_c $ \\
		\hline
		8292 & 16260 & \num{1.20E-3} & - & \num{1.21E-3} & - & \num{3.45E-2} & - \\
		\hline
		29069 & 57512 & \num{6.06E-4} & \num{9.88E-1} & \num{6.06E-4} & \num{9.98E-1} & \num{1.3E-2} & \num{1.41E+00} \\
		\hline
		115110 & 228970 & \num{3E-4} & \num{1.01E+00} & \num{3.01E-4} & \num{1.01E+00} & \num{4.6E-3} & \num{1.5E+00} \\
		\hline
		423918 & 845398 & \num{1.45E-4} & \num{1.05E+00} & \num{1.44E-4} & \num{1.06E+00} & \num{1.5E-3} & \num{1.62E+00} \\
		\hline
		1647132 & 3289420 & \num{6.88E-5} & \num{1.07E+00} & \num{6.86E-5} & \num{1.07E+00} & \num{4.56E-4} & \num{1.72E+00} \\
	\end{tabular}
\end{table}

\begin{table}[ht] 
	\caption{Test1. $L_{2}$ norms of errors for $d_{TL}$ refinement} \label{tab_test1_2}
	\centering
	\begin{tabular}{c|c|c|c|c|c}
		$N$ & $N_T$ & $\theta_{\epsilon}=\theta_K$ & $L_{2,u}$ & $L_{2,v}$ & $L_{\Psi}$ \\
		\hline
		115110 & 228970 & 800 & \num{3E-4} & \num{3.01E-4} & \num{4.6E-3}  \\
		\hline
		160753 & 329773 & 1600 & \num{1.98E-4} & \num{1.97E-4} & \num{3.01E-3}   \\
		\hline
		218507 & 435746 & 3200 & \num{1.178E-4} & \num{1.177E-4} & \num{1.7E-3}  \\
		\hline
		323158 & 645032 & 6400 & \num{6.01E-5} & \num{5.97E-5} & \num{7.6E-4}  \\
	\end{tabular}
\end{table}	

\subsection{Test 2. Comparison with analytical solution for Stokes flow regime in a three-layer channel} 

\begin{figure}[h!]
	\centering
	\begin{subfigure}{.4\textwidth}
		\centering
		\includegraphics[width=\textwidth]{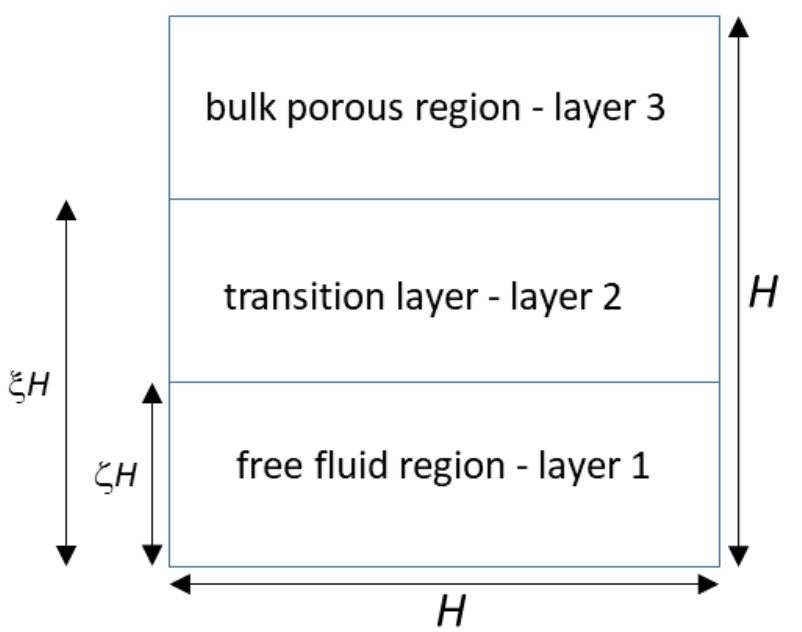}
		\caption{Setup}
		\label{fig_test2_1}%
	\end{subfigure}
	\begin{subfigure}{.5\textwidth} 
		\includegraphics[width=\textwidth]{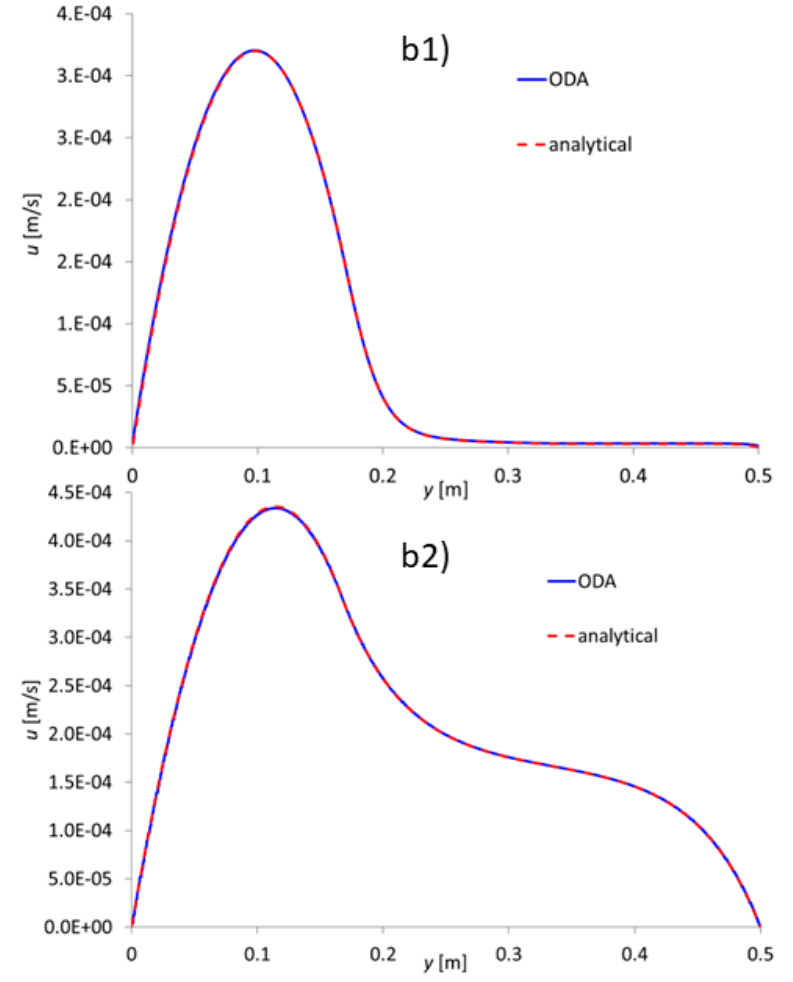}
		\caption{Numerical vs analytical solution. b1) ``ftl'' case with $Da=\num{2e-4}$. b2) ``ttl'' case with $Da=\num{1e-2}$}
		\label{fig_test2_2}%
	\end{subfigure}
	\caption{Test 2. Considered domain with three layers (a). Numerical and analytical results for two considered cases (b)}
	\label{fig_test2}%
\end{figure}

We deal with a 1D Stokes flow regime along the $x$ direction in a three-layer channel with total depth $H$, where the transition layer, TL, is in between a bulk porous and a clear fluid region (see \cref{fig_test2_1}). This problem is proposed in \cite{42}. The spatial distribution of the isotropic permeability $K$ is given in \cref{permb_test2}
\begin{subequations} \label{permb_test2}
	\begin{gather}
		\frac{1}{K\left(y\right)}=0 \qquad \textrm{in clear fluid region, layer 1},\qquad 0 \le y \le \zeta H,\\
		\frac{1}{K\left(y\right)}= \frac{1}{K_0} \frac{y-\zeta H}{\xi H - \zeta H} \qquad \textrm{in transition zone, layer 2} \qquad \zeta H \le y \le \xi H,  \\
		K\left(y\right)=K_0 \qquad \textrm{in bulk porous region, layer 3}, \qquad y > \xi H.
	\end{gather}
\end{subequations}	

The flow is driven by a uniform negative pressure gradient $G$ in the three layers, yielding the following Stokes-Brinkman Equations 
\begin{subequations} \label{test2_eqq}
	\begin{gather}
		\mu \frac{\partial^2 u}{\partial y^2}+G=0\qquad 0 \le y \le \zeta H,\\
		\mu \frac{\partial^2 u}{\partial y^2} - \frac{\mu u}{K\left(y\right)}+G=0\qquad \zeta H \le y \le \xi H,  \\
		\mu \frac{\partial^2 u}{\partial y^2} - \frac{\mu u}{K\left(y\right)}+G=0\qquad y > \xi H.
	\end{gather}
\end{subequations} 

The analytical solution of the dimensionless form of system in \cref{test2_eqq} is given by Eqq. (10) in \cite{42}. This has been obtained by solving the system in Eqq.(15)-(18) of the same paper. In the present work, we solved the system in Eqq. (15)-(18) in \cite{42} using the software package \textit{Mathematica} \cite{43}. We consider both \textit{thin} and \textit{fat} TL scenarios (``ttl'' and ``ftl'', respectively), and we present simulations for small and large Darcy numbers $Da$ values ($Da = K_0/H^2$). We set $H=0.5$ m and an equal value of the channel length. The ``ftl'' and ``ttl'' configurations are obtained by setting $\zeta = 1/3 H$ =  and $\xi= 2/3 H$ , and $\zeta= 1/3 H$  and $\xi= 1.006/3 H$, respectively. \par 
We set a pressure drop $\Delta \Psi = \num{5e-7} \textrm{ m}^2/\textrm{s}^2$ between the left and right end sides of the channel, while a no-slip velocity condition is assigned along the horizontal bottom and top walls. The kinematic fluid viscosity is $\nu=\num{1.5e-5} \textrm{ m}^2/\textrm{s}$. In \cref{tab4} we list the mesh sizes $\delta_0$ and $\delta_{TL}$ adopted for the different simulated scenarios, as well as the number of nodes $N$ and triangles $N_T$ of the corresponding meshes. The time step size $\Delta t$ = 10 s for all the simulated scenarios. We set $K_0$ such that $Da = \num{2e-4}$ and $Da = \num{1e-2}$, and $\textrm{CFL}_{max}$ ranges from 1.42 (``ttl'' case and small $Da$) to 2.45 (``ftl'' case and large Da). \par 
In \cref{fig_test2_2} we compare the numerical ODA solution against the analytical solution. The numerical solution of the presented solver  fits very well the analytical one. For brevity, we show only the solutions for ``ftl'' with small $Da$ and ``ttl'' with large $Da$ cases, but the matching is satisfactory also for the other investigated scenarios. The numerical solution in \cref{fig_test2_2} is related to mesh $m_1$ (see \cref{tab4}), but the solutions obtained for the meshes $m_0$ and $m_2$ are undistinguishable at the graphic scale, and for brevity are not shown. 

\begin{table} [h] \centering
	\caption{Test 2. Parameters for the computational meshes $m_i$}  \label{tab4}
	\begin{tabular}{c|c|c|c|c}
		& $\delta_0$ [m] & $\delta_{TL}$ [m] & N & $N_T$ \\
		\hline
		ftl $m_0$ & \num{5e-3} & \num{5e-3} & 19478 & 37752 \\
		\hline
		ftl $m_1$ & \num{2.5e-3} & \num{2.5e-3} & 51841 & 102287 \\
		\hline
		ftl $m_2$ & \num{1.25e-3} & \num{1.25e-3} & 195556 & 389356 \\
		\hline
		ttl $m_0$ & \num{5e-3} & \num{7e-5} & 242750 & 484234 \\
		\hline
		ttl $m_1$ & \num{2.5e-3} & \num{7e-5} & 273272 & 545112 \\
		\hline
		ttl $m_2$ & \num{1.25e-3} & \num{7e-5} & 402484 & 803165 \\
	\end{tabular}
\end{table}	

\subsection{Test 3. Comparison with pore-scale results for Stokes flow regime}  
\begin{figure}[h!] \label{fig_test3_1}
	\centering
	\begin{subfigure}{.55\textwidth}
		\includegraphics[width=\textwidth]{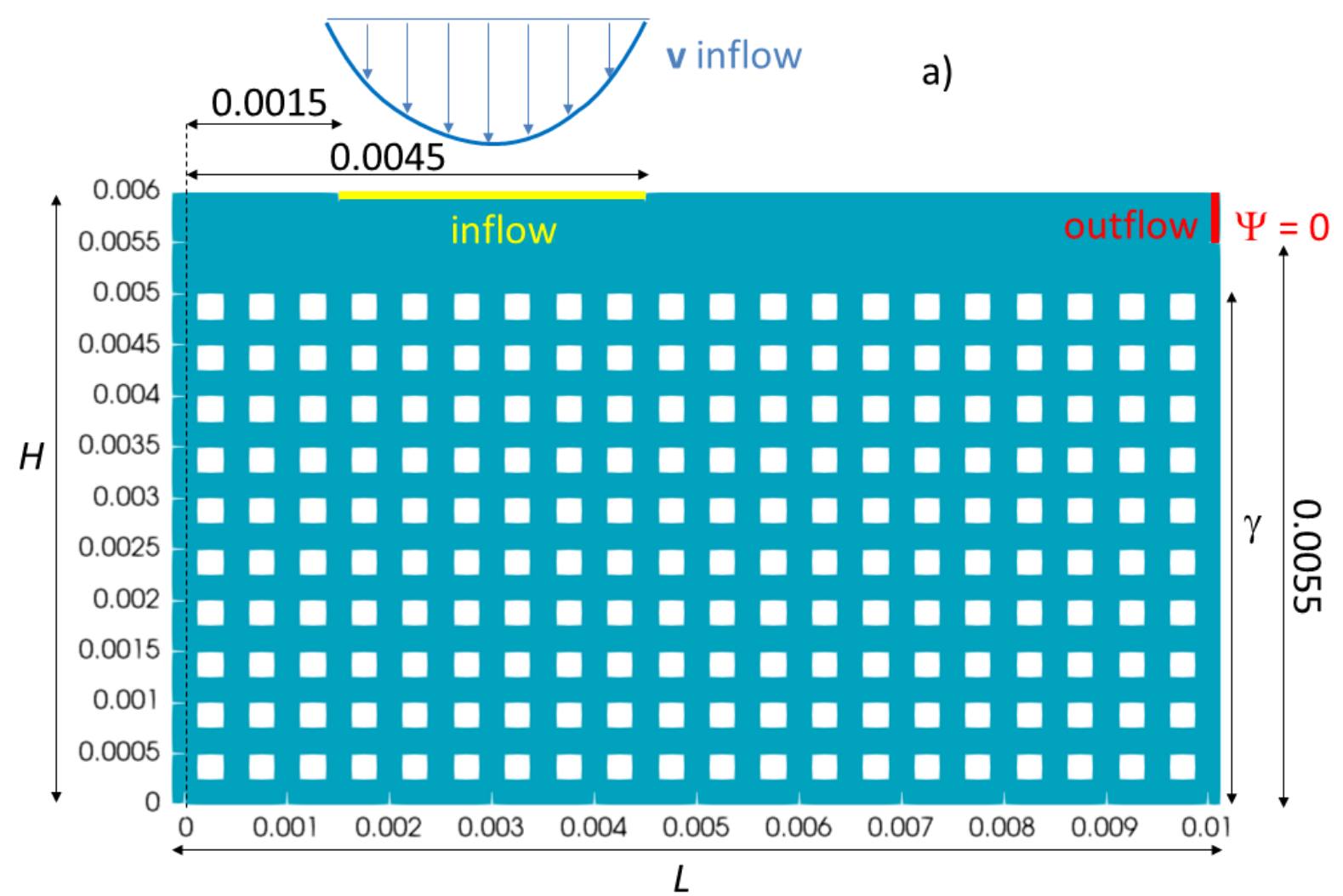}
		\caption{Geometrical pore-scale setup and boundary condition (measures in m)}
		\label{fig_test3_1_a}%
	\end{subfigure}
	\hfill 
	\begin{subfigure}{.4\textwidth} 
		\includegraphics[width=\textwidth]{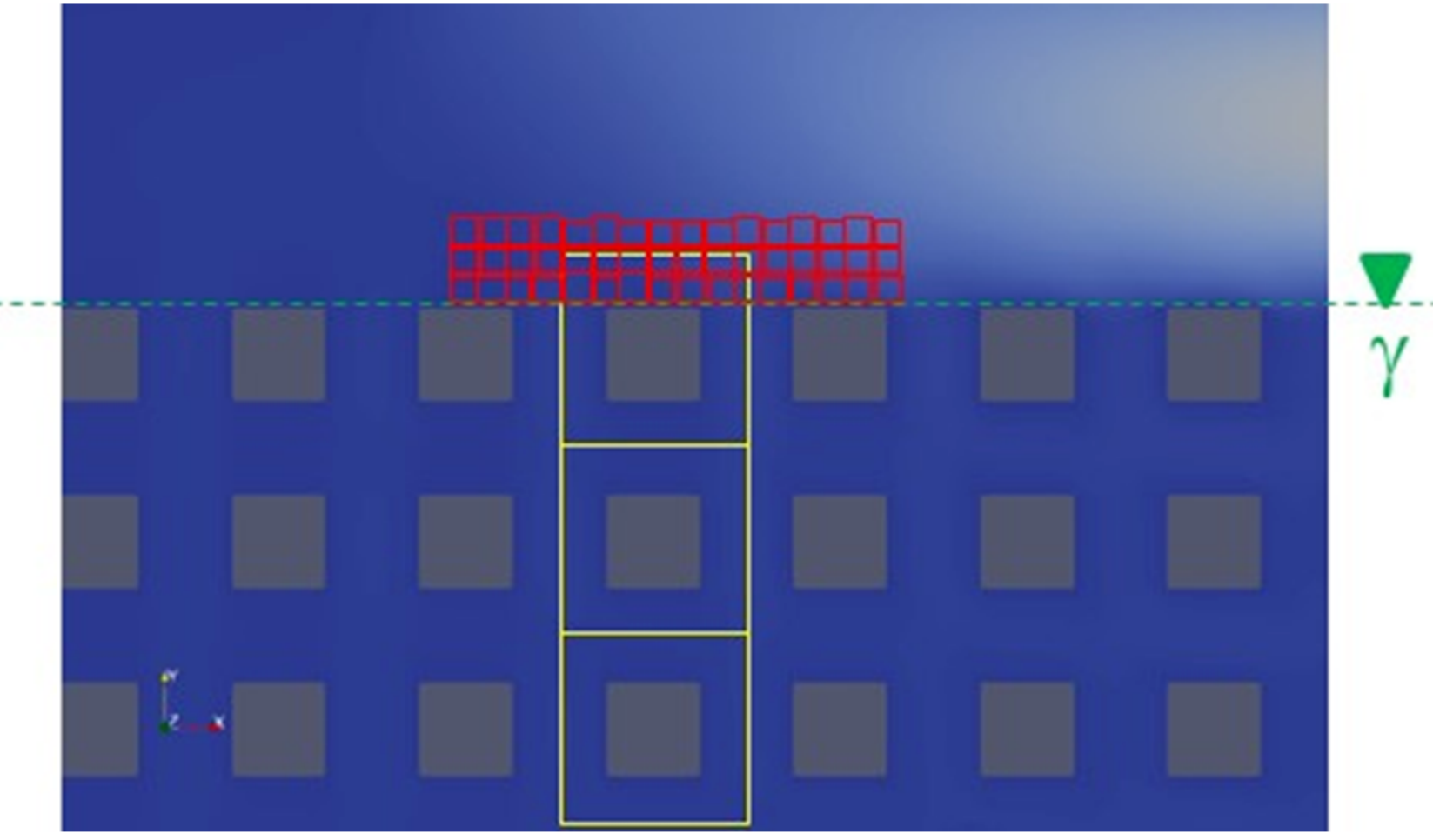}
		\caption{REVs within $\Omega_{pm}$ and $\Omega_{ff}$}
		\label{fig_test3_1_b}%
	\end{subfigure}
	\caption{Test 3. Pore-scale setup (a). Considered REVs for averaging (b)}		
\end{figure}

Test 3 is related to filtration/exfiltration processes in microfluidics, which presents several industrial, environmental or biomedical applications (e.g., membrane filtration processes, wastewater treatment, water purification, hemodialysis, ... ). The test case considered in this Section has been proposed in \cite{44}. The computational domain $\Omega=(0, L) \times (0, H)$ has a free fluid region $\Omega_{ff} = (0, L) \times (\gamma,H)$ and a porous region $\Omega_{pm} = (0,L)\times (0,\gamma)$ with $L=10.25$ mm and $H = 6$ mm. The porous region is isotropic, made up of 20 $\times$ 10 square solid inclusions of size $d_{si}= \SI{250} {\micro\metre}$, such that the porosity $\epsilon$ is 0.75. $\gamma$ = 5 mm is the distance between the bottom boundary and the tangential line located on top of the uppermost row of solid inclusions. This is shown in \cref{fig_test3_1_a}, where the geometrical setup and the assigned boundary conditions are depicted. An inflow velocity profile $\mathbf{v}=(0,V_{max}sin(\frac{1000}{3} \pi x))$ is assigned on $\Gamma_e = 1.5 (\textrm{ mm}, 4.5 \textrm{ mm}) \times H$, with $V_{max}=\num{1e-3}$ m/s, and $\Psi = 0$ on $\Gamma_n = L \times (5.5 \textrm{ mm}, H)$. The remaining portions of the boundaries are assumed to be impervious walls with no-slip velocity condition. The considered fluid is water ($\nu = \num{1e-6} \textrm{ m}^2/\textrm{s}$) and the bulk  permeability is $K_0= \num{3.45e-9} \textrm{ m}^2$, computed according to the geometry of the solid inclusions \cite{45}. The Reynolds number associated with the free fluid zone, $Re = \left(V_{max} \left(H-\gamma\right)\right) / \nu$, is equal to 1.\par
To obtain a reference solution, we averaged the results of a pore-scale simulation (PSS) obtained by solving the present test case using the open-source simulator $\textrm{DuMu}^\textrm{x}$ \cite{40}. At the pore scale, the flow is governed by the Stokes equations in the entire domain, with no-slip BCs assigned on the boundaries of the solid inclusions. A uniform structured mesh was used for the PSS, with side $\delta_{PSS}=\num{1.5625e-5}$ m. The ODA solution is compared with the ``surface average'' \cite{11, 12} results of the PSS over square-shaped REVs, the size of which, $l_{REV}$, depend on whether the REV is located in the porous or the fluid region (see \cref{fig_test3_1_b}). Within $\Omega_{pm}$, $l_{REV}= 2 d_{si}$ and REVs centroids coincide with those of the solid inclusions (yellow squares in \cref{fig_test3_1_b}). Within $\Omega_{ff}$, $l_{REV}= 8 \delta_{PSS}$ (red squares in \cref{fig_test3_1_b}). \par 
Both \textit{bottom} and \textit{middle} TL positions were assumed. \textit{Bottom} and \textit{middle} positions mean that the center-line of the TL is at $y=\gamma-d_{si} / 2$ or $y=\gamma$, respectively. The spatial variations of $\epsilon$ and $1/K$ within TL are given by \cref{eq:3}. \par
The TL width, $d_{TL}$, as well as the parameters needed for modelling porosity and permeability variation given by \cref{eq:3} were selected according to preliminary simulations performed on a fine mesh (uniform linear size of 0.0085 mm, $N$=895064 and $N_T$ = 1786418), and were adopted for TL positions. Therefore we assumed: $d_{TL} \in \left[d_{si}, 3 d_{si}\right]$, $d_{\epsilon} = 0$, $d_K \in \left[0,\num{1e-4} \textrm{m}\right]$ and $\theta_{\epsilon} = \theta_K \in \left[ \num{7e+3}, \num{5e+4}\right]$. For each simulation we compared, along the horizontal line at $y = 0.004875$ m (the vertical coordinate of the centers of mass of the solid inclusions of the uppermost row), the distributions of $u$, $v$ and $\Psi$ computed by the ODA solver with the corresponding reference data of the averaged PSS.\par
The ODA solutions obtained for the \textit{middle} TL position seemed not to be sensitive to the size $\delta_{TL}$, and, compared to the \textit{bottom} position configuration, we observed 1) an overestimation of the velocity vector magnitude $ \lVert \textbf{u} \rVert$, up to twice in the bulk region$\Omega_{pm}$ and up to 20 \% in the $\Omega_{ff}$ region, and 2) an underestimation of $ \lVert \textbf{u} \rVert$ close to the interface between $\Omega_{pm}$ and $\Omega_{ff}$, up to one magnitude order. This is why the results are not further discussed here. The best match with the averaged PSS results was obtained for thr \textit{middle} TL configuration with $d_K = \num{3e-5}$ and $\theta_{\epsilon} =\theta_K = \num{4.7e+4}$.\par 	
The adopted mesh sizes are $\delta_0=\num{2.5e-5}$ m and $\delta_{TL}=\num{1e-5}$ m. These sizes guarantee a good compromise between the accuracy of the results (maximum relative error values of the velocity components and kinematic pressure not greater than 1 \% compared to the results obtained over the fine mesh above mentioned) and the computational effort. The time step size is $\num{2.5e-3}$ s and $CFL_{max} \simeq 1.98$.\par
\begin{figure}[h!]
	\centering
	\includegraphics[width=1\textwidth]{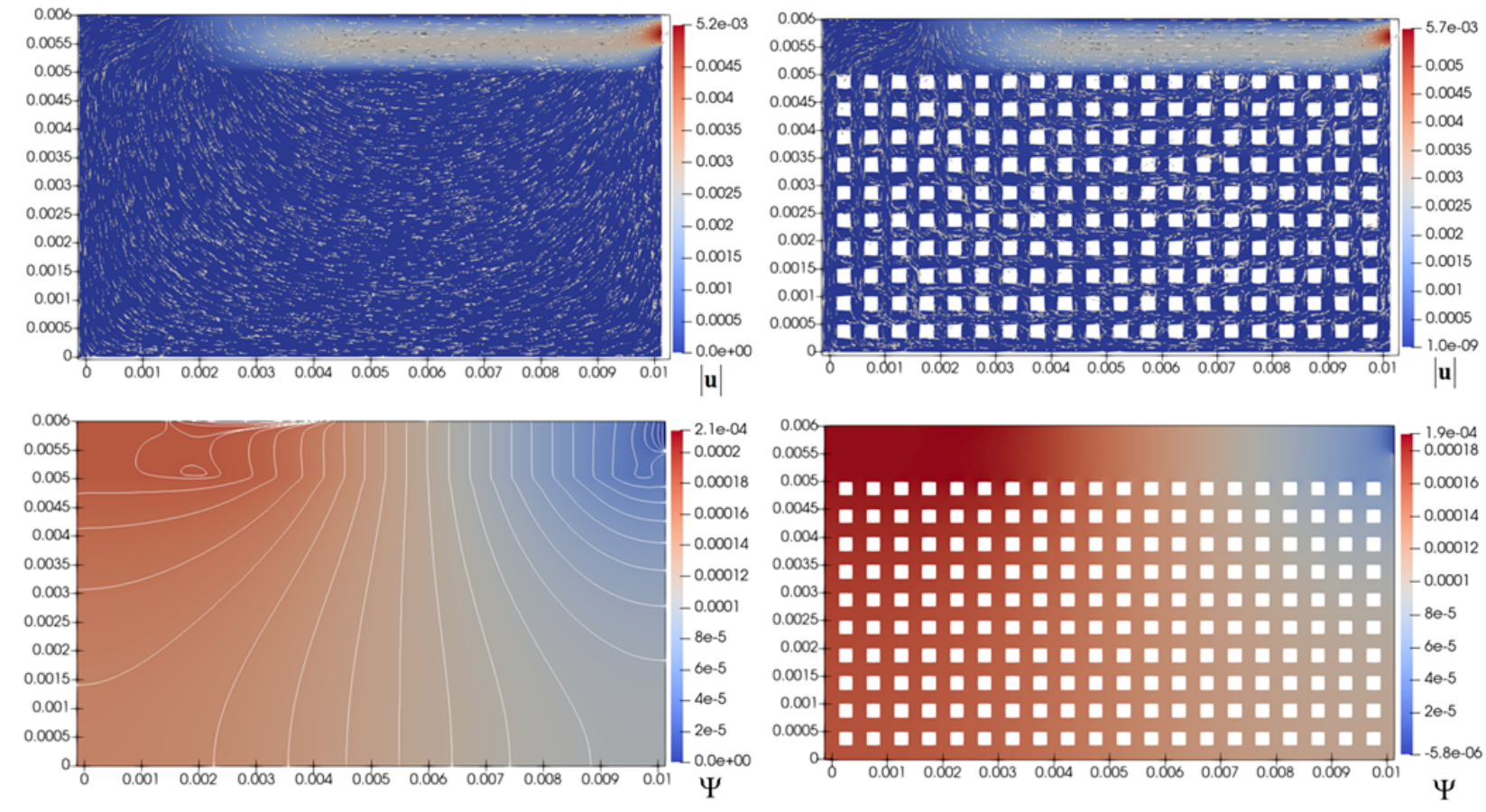}
	\caption{Test3. Comparison between the velocity and kinematic pressure fields computed by the ODA solver (left column) and the PSS results \cite{40} (right column). Top row: velocity (in m/s), bottom row: kinematic pressure (in $\textrm{m}^2/\textrm{s}^2$)}
	\label{fig_test3_2}
\end{figure}

In \cref{fig_test3_2} we compare the velocity and kinematic pressure fields provided by the present ODA solver and the corresponding PSS results. The porous medium is an obstacle to the incoming flow from $\Gamma_e$, which largely deviates in the upper free fluid region, to the right, forming a channelized flow. Part of the incoming flow infiltrates the left portion of $\Omega_{pm}$ and exfiltrates through the right side, towards the outlet of the domain. The absolute value of the velocity vector in the bulk $\Omega_{pm}$ is approximately 1.5-2 orders of magnitude  smaller than in $\Omega_{ff}$. The ODA solver well reproduces the overall flow and pressure fields predicted by the PSS, with a small overestimation of $\Psi$ close to the inflow region. \par
\begin{figure}[h!]
	\centering
	\includegraphics[width=1\textwidth]{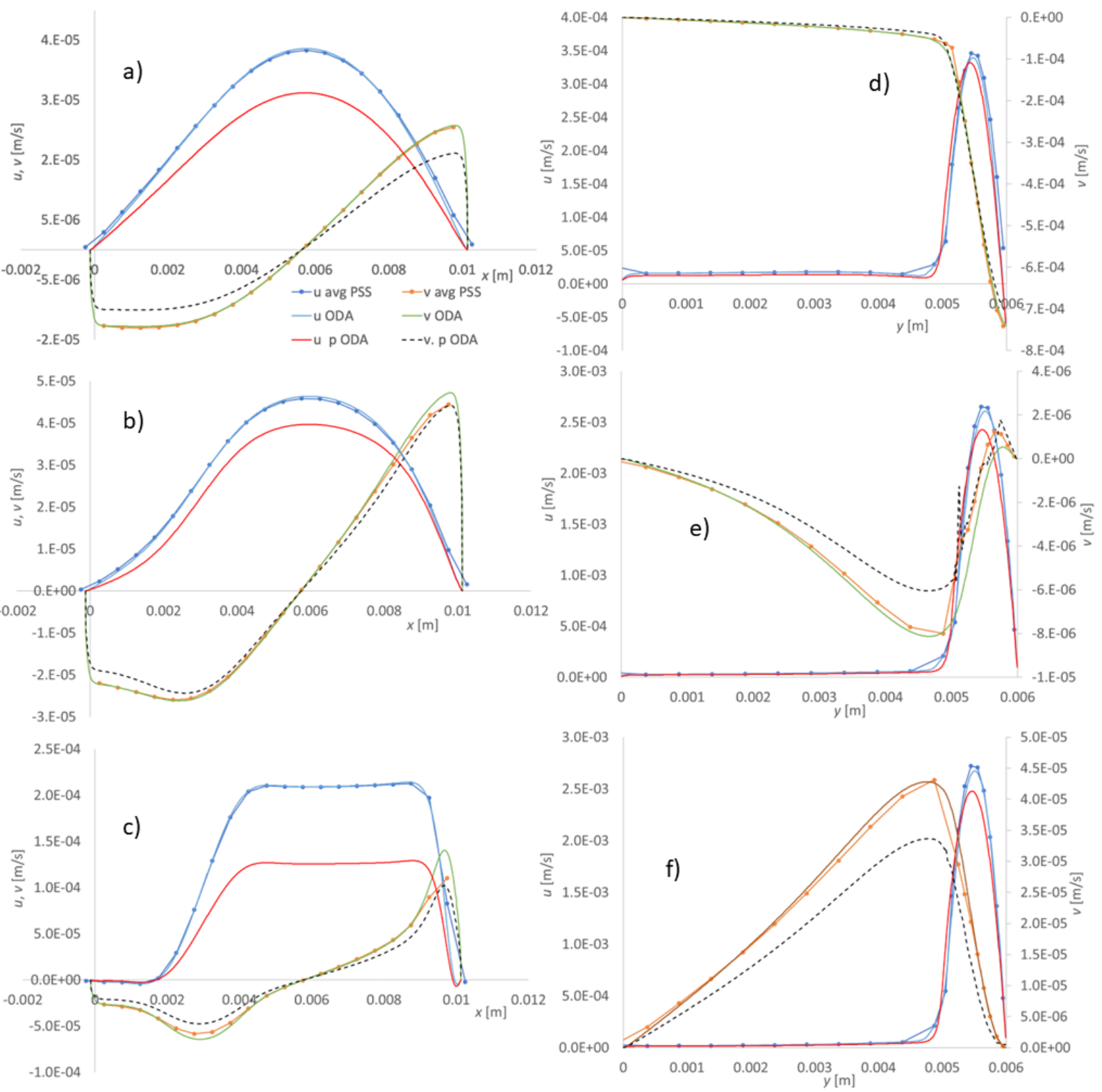}
	\caption{Test3. Velocity components computed by the averaged PSS, the ODA solver and the penalized ODA solver (p ODA) at different positions. a) $y = 0.001875$ m; b) $y = 0.003375$ m; c) $y = 0.004875$ m; d) $y = 0.00225$ m; e) $y = 0.00525$ m; f) $y = 0.00825$ m}
	\label{fig_test3_3}
\end{figure}
\begin{figure}[h!]
	\centering
	\includegraphics[width=1\textwidth]{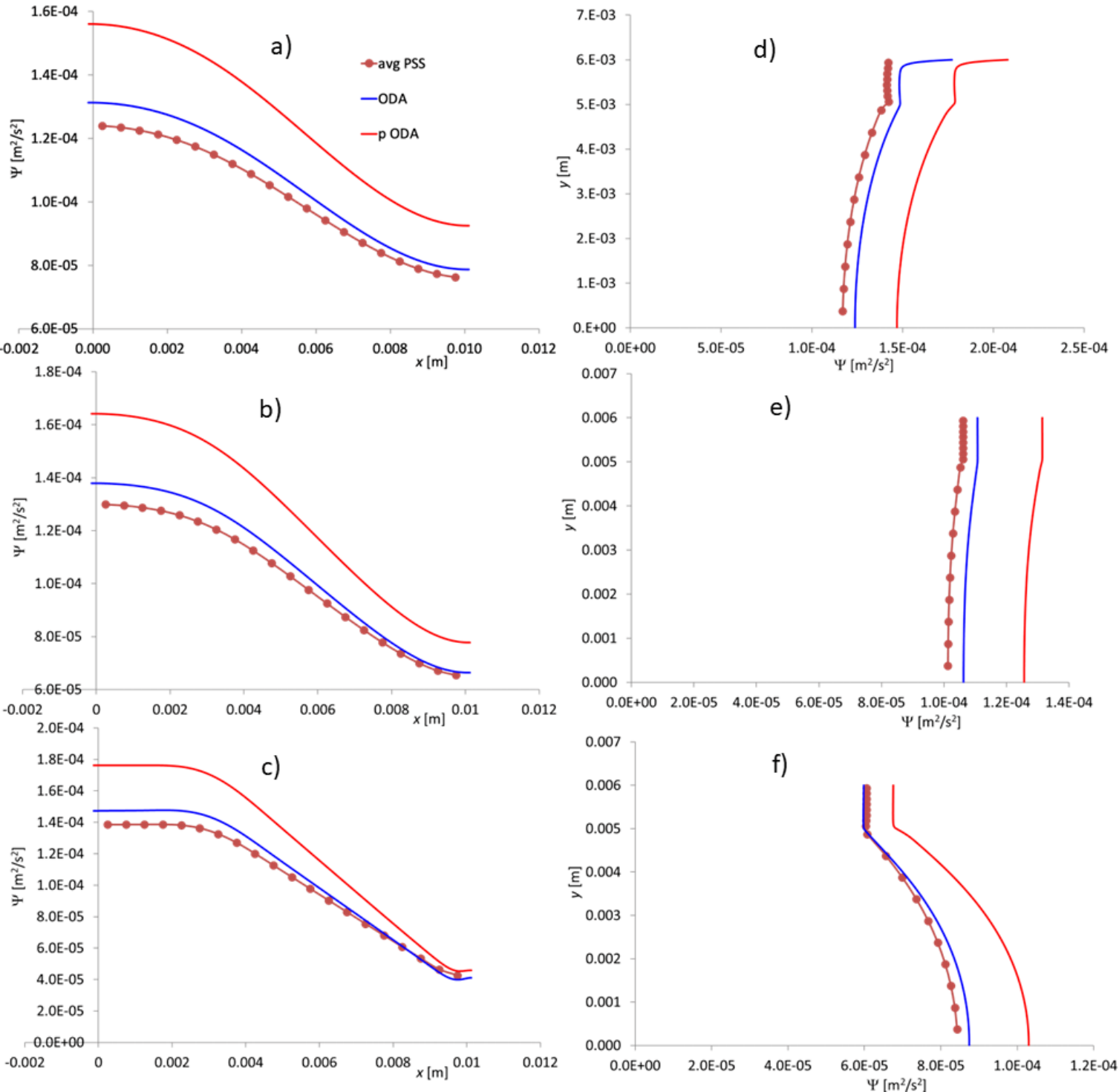}
	\caption{Test3. Kinematic pressure computed by the averaged PSS, the ODA solver and the penalized ODA solver (p ODA) at different positions. a) $y = 0.001875$ m; b) $y = 0.003375$ m; c) $y = 0.004875$ m; d) $y = 0.00225$ m; e) $y = 0.00525$ m; f) $y = 0.00825$ m}
	\label{fig_test3_4}
\end{figure}

Velocity and kinematic pressure profiles are compared in \cref{fig_test3_3} and \cref{fig_test3_4}. There, excellent agreement between the ODA and the averaged PSS solutions is observed within $\Omega_{pm}$ ($y=0.001875$ m  and $y=0.003375$ m) and close to the interface ($y=0.004875$ m). According to the computed $v$ component and $\Psi$ at $y=0.001875$ m, $y=0.003375$ m and $y=0.004875$ m, as well as at $x=0.00225$ m and $x=0.00825$, the proposed solver correctly predicts the infiltration processes into and exfiltration processes from the porous domain. The ODA solver slightly underestimates the peak value of the $u$ component in the $\Omega_{ff}$ region compared to the averaged PSS (positions $x = 0.00225$ m, $x = 0.00525$ m and $x = 0.00825$ m).  \par 
In \cref{fig_test3_3} and \cref{fig_test3_4}, we also plot the solution of a \textit{penalized} ODA solver, (see \cref{intro} and \cite{17, 18, 19, 20, 21, 22} for further details). This solver is obtained by assuming a discontinuous function for $\epsilon$ and $1/K$ at the transition between $\Omega_{pm}$ and $\Omega_{ff}$, i.e., at $y=\gamma$ (the associated results are marked with ``p ODA''). This strongly affects the infiltration/exfiltration processes, as observed in these figures. The absolute value of the $v$ component is underestimated close to the interface and along the three vertical profiles at $x=0.00225$ m, $x=0.00525$ m and $x=0.00825$ m. Furthermore, unphysical oscillations in the $v$ profile at $x=0.00525$ m are observed due to the interfacial stress jump. Significant underestimation of the $u$ component within the porous medium, in the free fluid region and close to the interface, are also observed. The kinematic pressure is overestimated throughout the computational domain. This analysis shows that the ``penalized'' approach provides a poor estimation of the results, not only close to the interface but also within the bulk $\Omega_{pm}$ and $\Omega_{ff}$ regions, since the overall infiltration/exfiltration processes are not properly recovered, due to the choice of the discontinuous profiles of $\epsilon$ and $1/K$. 

\subsection{Test 4. Comparison with pore-scale results for Navier-Stokes flow regime}  
\begin{figure}[h!]
	\centering
	\includegraphics[width=0.4\textwidth]{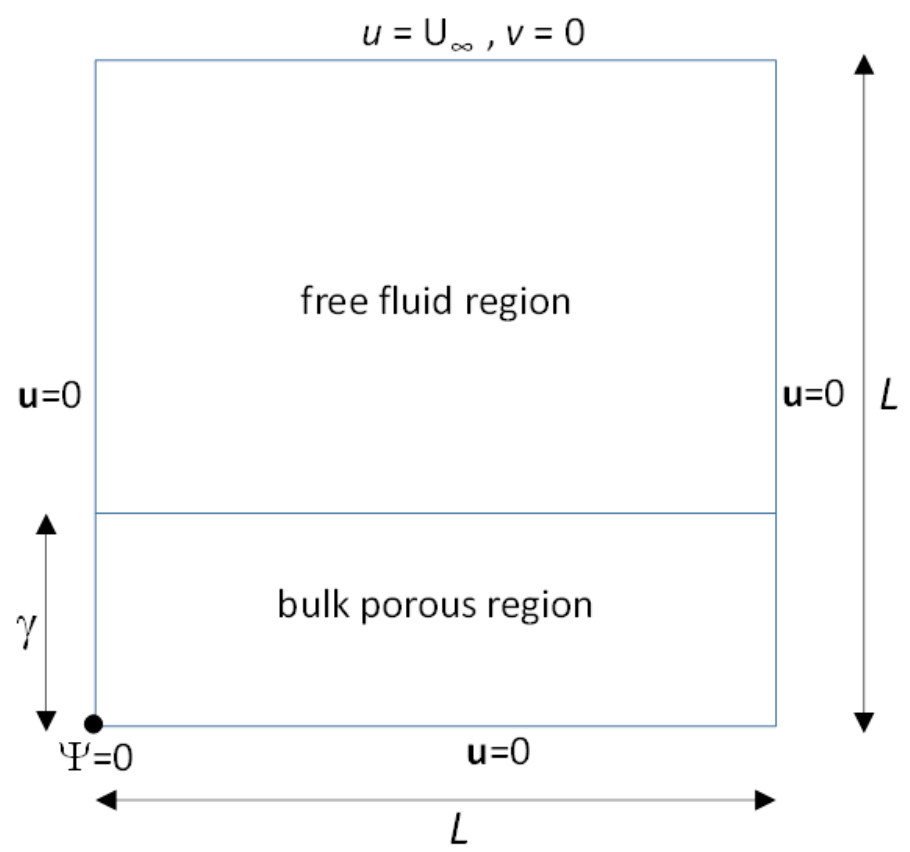}
	\caption{Test4. Computational domain and assigned boundary conditions, with $L$ = 1 m and $\gamma = L/3$}
	\label{fig_test4_1}
\end{figure}

In test 4, we deal with a lid-driven cavity flow over a fibrous porous medium. The lid-driven cavity flow is an idealized paradigm of internal flows of industrial or natural processes, for example industrial microelectronics, metal casting, flows over slots on the walls of heat exchangers, dynamics of lakes, ... The vortical flow structure and the related momentum transport process can be modulated by porous medium, which can be used as a passive flow control tool. 
The present test case has been proposed in \cite{6}. Here we deal with a square domain $\Omega=(0, L)^2$, consisting of an upper free fluid region $\Omega_{ff} = (0, L) \times (\gamma,L)$ and a bottom porous region $\Omega_{pm} = (0,L)\times (0,\gamma)$, made of vertical fibers, with porosity $\epsilon$ = 0.8, and orthotropic tensor. The setup is shown in \cref{fig_test4_1}, with $L=1$ m and $\gamma = L/3$. The upper boundary moves to the right with with assigned horizontal velocity $U_{\infty}=1 \textrm{ m/s}$; the other boundaries are impervious. We ran simulations for $Re$ = 100 and 1000, where $Re=\frac{U_\infty L}{\nu}$. The associated pore-scale Reynolds number within $\Omega_{pm}$ ranged between \num{1.75e-2} and \num{5.125e-3} \cite{6}. \par
The authors of \cite{6} performed Direct Numerical Simulations (DNS) at the pore scale, where the Navier Stokes equations with no-slip BCs on the fibers boundaries are solved, and the results have been then averaged over cubic REV volumes with side length $l_{REV}=0.02$ m. From the averaged DNS results, they also computed the bulk permeability tensor coefficients $K_{i,j}^0, i,j = 1,2$, as well as the extensions $d_{TL_i}$, measured from the top of $\Omega_{pm}$, where $1\big/ K_{i,j}$ decreases, almost linearly, from the bulk value to zero. More details can be found in the aforementioned paper. In \cref{tab_test_4.1} we list the values of $K_{i,j}^0$ and $d_{TL_i}$. Due to the orthotropy of the porous region, $K_{i,j}^0 = 0$ with $ i\neq j$. The averaged DNS simulations are the reference solutions for the comparison with the proposed ODA solver.\par 
We simulated different scenarios, where the position and the size of the transition layer $d_{TL}$ is changed (see \cref{tab_test_4.2}). As before, with the nomenclature \textit{bottom} and \textit{middle} we refer to the position of the TL, whose top level is at $y=\gamma$ and $y=\gamma+d_{TL}/2$, respectively. We assumed linear variation of the porosity along $d_{TL}$, and if $d_{TL_1} \neq d_{TL_2}$, we assumed linear variation of $\epsilon$ along the average distance $\left(d_{TL_1} + d_{TL_2}\right)/2$. We performed a sensitivity analysis to the mesh size, where the reference solution was obtained over a fine mesh (uniform mesh size 0.001 m in the entire domain, $N_T$=2241989, $N$ = 1122965). The sizes $\delta_0$ and $\delta_{TL}$ m listed in \cref{tab_test_4.2} provided a good compromise between the accuracy of the results (maximum value of the relative errors of the velocity components and kinematic pressure compared to the results of the fine mesh smaller than 1 \%) and the computational costs. The time step size of the numerical simulation is $\Delta t$ = 0.02 s and $CFL_{max}$ ranges from 3.15 ($Re = 1000$) to 3.24 ($Re = 100$). \par
\begin{figure}[h!]
	\centering
	\includegraphics[width=0.8\textwidth]{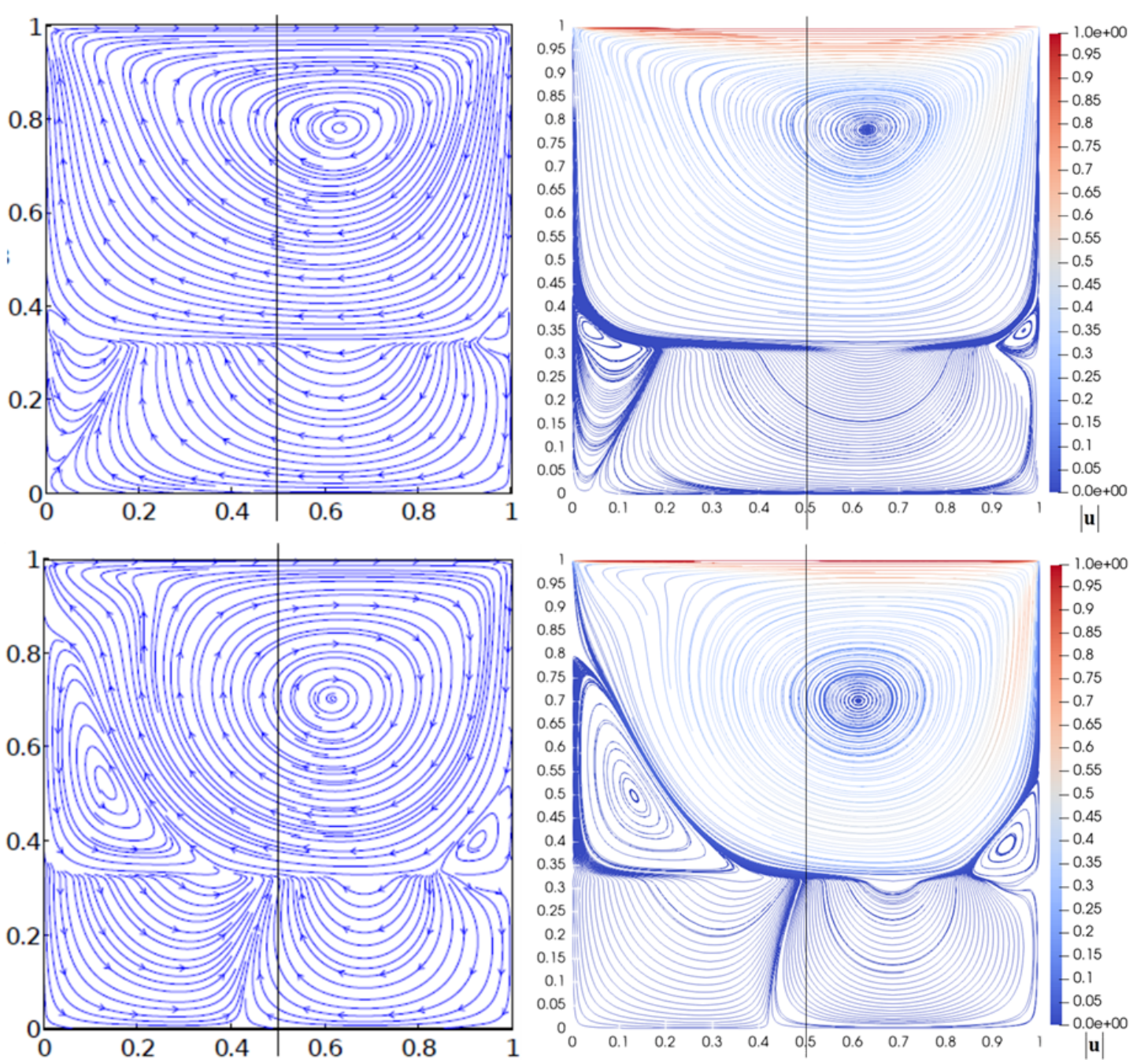}
	\caption{Test4. Comparison between the velocity streamlines computed by the ODA solver and the averaged DNS results \cite{6}. Left column: averaged DNS results, right column: ODA results. Top row: $Re = 100$, bottom row: $Re = 1000$ (velocity of the proposed model in m/s)}
	\label{fig_test4_2}
\end{figure}
\begin{figure}[h!]
	\centering
	\includegraphics[width=0.8\textwidth]{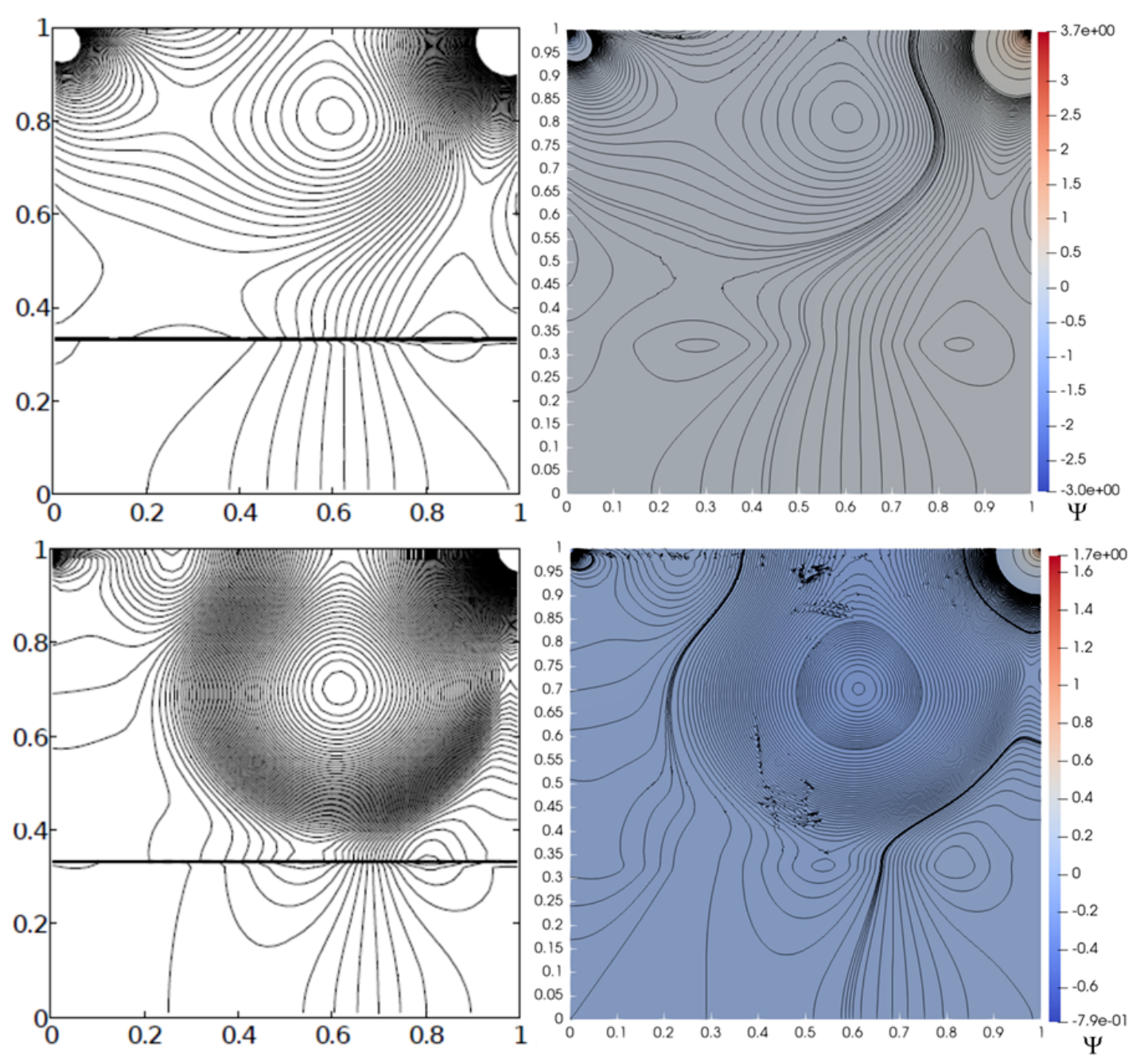}
	\caption{Test4. Comparison between the contours of the (kinematic) pressure computed by the ODA solver and the averaged DNS results \cite{6}. Left column: averaged DNS results, right column: ODA results. Top row: $Re = 100$, bottom row: $Re = 1000$ (kinematic pressure of the proposed model in $\textrm{m}^2/\textrm{s}^2$)}
	\label{fig_test4_3}
\end{figure}

In \cref{fig_test4_2} and \cref{fig_test4_3} we compare the velocity streamlines and the pressure contours computed by the presented ODA solver with the averaged reference DNS results of \cite{6}. Excellent agreement can be found for both $Re$ values. As $Re$ increases, the large vortex within $\Omega_{ff}$ moves downwards to the center of the domain, and the size of the corner vortices increases. According to the streamlines, we argue that the porous region represents an obstacle for the flow, and only a small amount of the fluid penetrates $\Omega_{pm}$, with a velocity magnitude approximately three orders of magnitude smaller than in $\Omega_{ff}$. The separation between the two major recirculation zones within $\Omega_{pm}$ moves to the right as $Re$ increases from 100 to 1000. The local pressure minima are associated with the centers of the vortices, both within $\Omega_{ff}$ and $\Omega_{pm}$. \par
\begin{figure}[h!]
	\centering
	\includegraphics[width=1\textwidth]{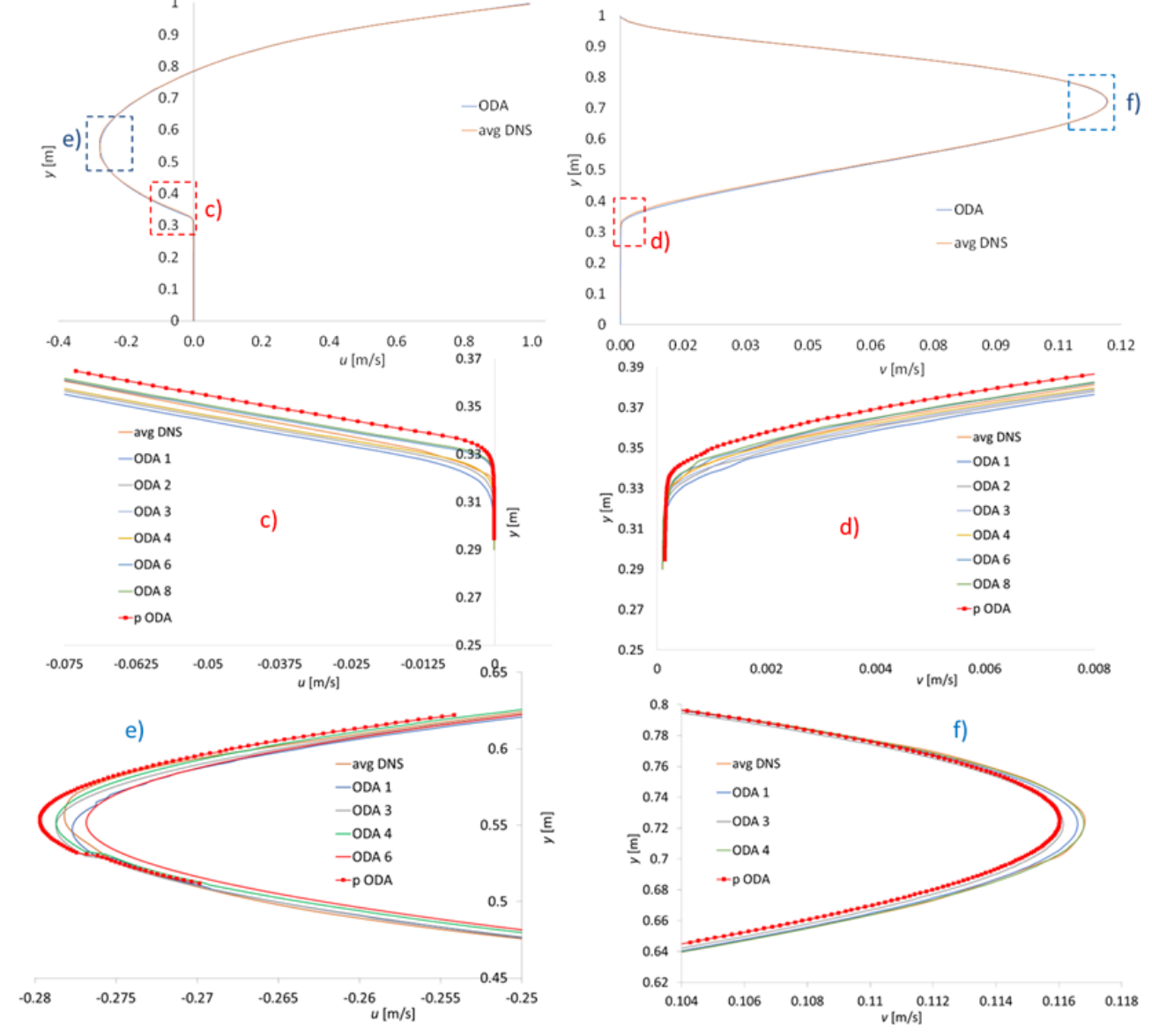}
	\caption{Test4. Velocity components computed by the ODA solver and the averaged DNS results \cite{6} at $x$=0.5 m, $Re$ = 100. a) $u$; b) $v$ ; c) $u$, zoom in the peak value region;  d) $v$, zoom in the peak value region; e) $u$, zoom in the transition layer region; f) $v$, zoom in the transition layer region}
	\label{fig_test4_4}
\end{figure}

\begin{figure}[h!]
	\centering
	\includegraphics[width=1\textwidth]{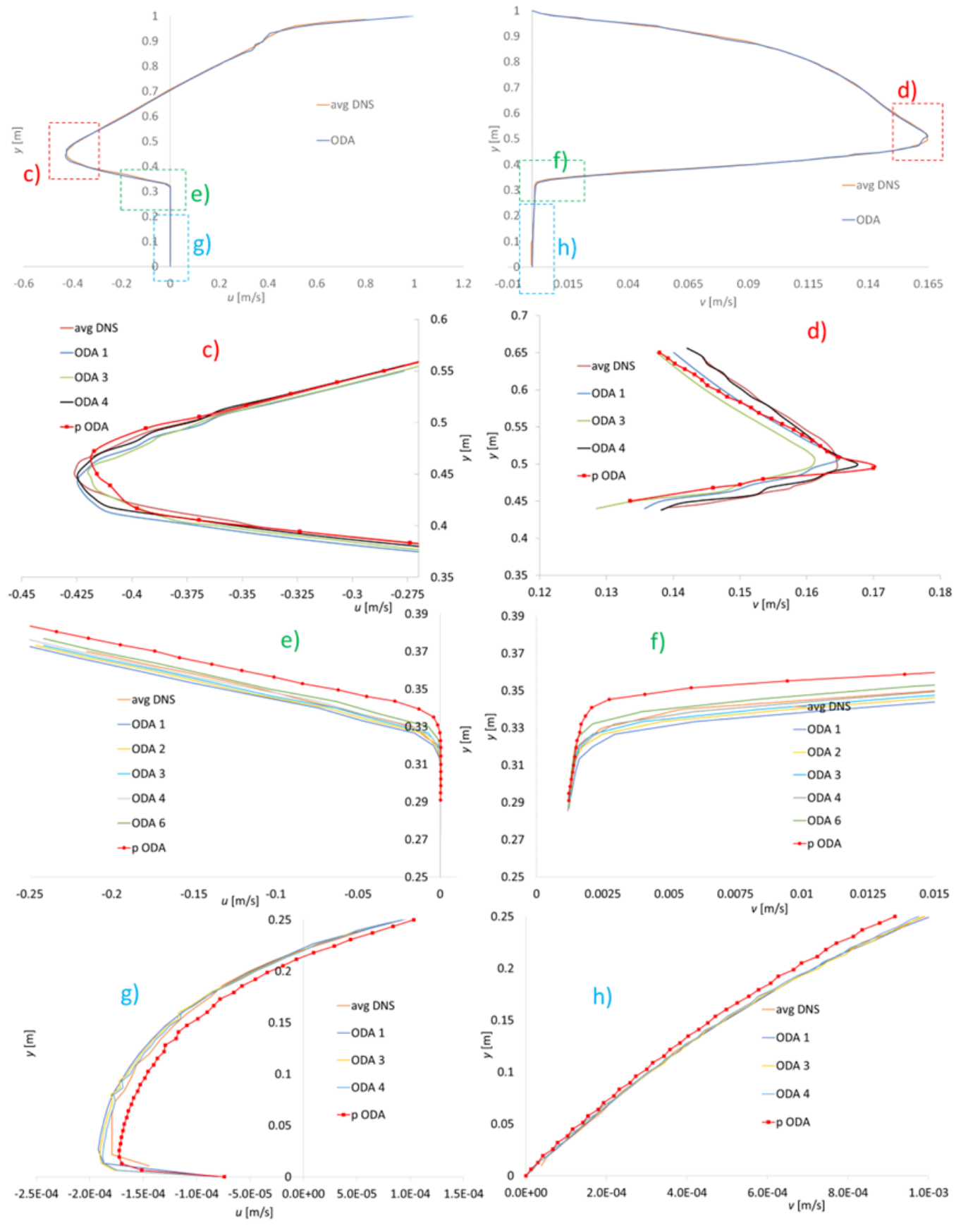}
	\caption{Test4. Velocity components computed by the ODA solver and the averaged DNS results \cite{6} at $x$=0.5 m, $Re$ = 1000. a) $u$; b) $v$ ; c) $u$, zoom in the peak value region;  d) $v$, zoom in the peak value region; e) $u$, zoom in the transition layer region; f) $v$, zoom in the transition layer region; g) $u$, zoom in the bulk $\Omega_{pm}$ region; h) $v$, zoom in the bulk $\Omega_{pm}$ region}
	\label{fig_test4_5}
\end{figure}

 
Overall good agreement between the present model and the reference solution can be observed for the velocity components computed along the vertical center-vertical-line of the domain ($x$ = 0.5 m) (see \cref{fig_test4_4} and \cref{fig_test4_5}). The most accurate results are obtained by setting a middle position of the transition layer and the distance $d_{TL_1}$ and $d_{TL_2}$ to be equal to the characteristic size of the REV used for the averaging process of the DNS results (scenario 4 in \cref{tab_test_4.2}). In the same \cref{fig_test4_4} and \cref{fig_test4_5}, the solutions marked as ``p ODA'' are obtained in the framework of a \textit{penalized} approach, with a discontinuity of $\epsilon$ and $1/K_{i,i}, i=1,2$ across the interface placed at $y=\gamma$. Due to the interfacial stress jump, $u$ and $v$ profiles close to interface between $\Omega_{pm}$ and $\Omega_{ff}$ are shifted. Poor estimation of the peak values is also observed in $\Omega_{ff}$, and both velocity components are underestimated within $\Omega_{pm}$. 

\begin{table} [h!]
	\caption{Test 4. Permeability coefficients (from \cite{6})} \label{tab_test_4.1}
	\centering
	\begin{tabular}{c|c|c|c|c}
		& $K_{1,1}^0$ [m$^{2}$]  &$K_{2,2}^0$ [m$^{2}$] & $d_{TL_1}$ [m]   & $d_{TL_2}$ [m]\\
		\hline
		$Re$ 100 & \num{1.05e-5} & \num{2.196e-2} & 0.05  & 0.03 \\
		\hline
		$Re$ 1000 & \num{1.06e-5} & \num{2.25e-5} & 0.05   & 0.03 \\
	\end{tabular}
\end{table}

\begin{table} [h!]
	\caption{Test 4. Mesh parameters for different scenarios} \label{tab_test_4.2}
	\centering
	\begin{tabular}{c|c|c|c|c|c|c}
		scenario  & position TL  & $d_{TL_1}$ [m] & $d_{TL_2}$ [m] & $\delta_{TL}$ [m] & $N$ & $N_T$  \\
		\hline
		1 & bottom & 2.5 \textit{l} & 1.5 \textit{l} & 0.0015  & 241187 & 480781 \\
		\hline
		2 & bottom & 1.5 \textit{l}  & 1.5 \textit{l} & 0.0015 & 235282 & 468978 \\
		\hline
		3 & bottom &  \textit{l}  & \textit{l} & 0.0015 & 233128 & 464673 \\
		\hline
		4 & middle & \textit{l}  & \textit{l}  & 0.0015 & 233948 & 466351 \\
		\hline
		5 & bottom & 0.5 \textit{l}  & 0.5 \textit{l}  & 0.00075 & 210725 & 419851 \\
		\hline
		6 & middle & 0.5 \textit{l} & 0.5 \textit{l} & 0.00075 & 210839 & 419941 \\
		\hline
		7 & bottom & 0.1 \textit{l} & 0.1 \textit{l} & 0.00015  & 359986 & 718352 \\
		\hline
		8 & bottom & 0.05 \textit{l} & 0.05 \textit{l} & 0.000075  & 496636 & 991647 \\
		\hline
		9 & bottom & 0.025  \textit{l} & 0.025 \textit{l} & 0.000035 & 793656 & 1585690 \\
	\end{tabular}
\end{table}

\subsection{Test 5. Analysis of free fluid flow over an anisotropic porous obstacle for different $Re$ values}   	
\begin{figure} [h!]
	\centering
	\includegraphics[width=0.67\textwidth]{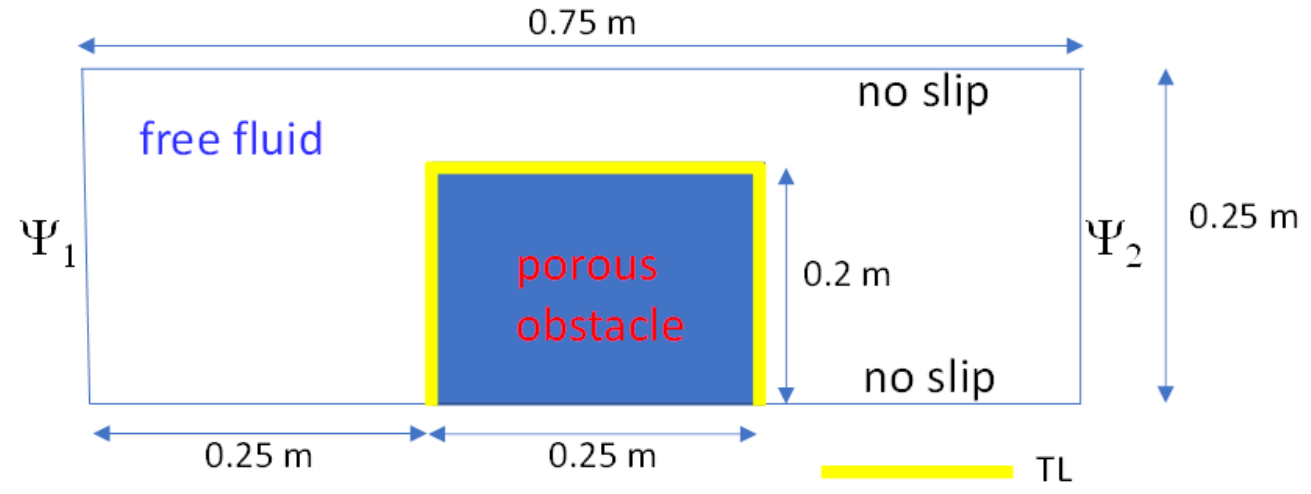}
	\caption{Test 5. Definition sketch ($Re \ll 1$ case) and assigned BCs}
	\label{figure_5_1}
\end{figure}

A porous obstacle invested by a fluid finds several applications (e.g., oil filters, porous coating acting as passive flow control device, porous regulating flow devices, ...). The test case proposed in this Section has been proposed in \cite{46}. Here a free fluid flows around and within an anisotropic porous obstacle $\Omega_{pm}$ with bulk porosity 0.4 (see \cref{figure_5_1}). Flow is driven by a pressure drop $\Delta \Psi = \Psi_1- \Psi_2$ between the upstream and downstream sides of the domain, while no-slip velocity BC is imposed on the bottom and top sides. This test is proposed in \cite{46}. The kinematic viscosity of the fluid is $\nu=\num{1.5e-5} \textrm{ m}^2/s$. The Reynolds number is calculated as $Re=\frac{\| \mathfrak{u}\|_{max}h_{ff}} {\nu}$, where $\| \mathfrak{u} \|_{max}$ is the maximum value of the velocity vector magnitude in the free fluid region above the porous obstacle and $h_{ff}$ is the fluid depth above the obstacle. We simulate two cases of $Re\ll 1$ and $Re \simeq 160$. \par  
The anisotropic tensor $\mathbf{K}$ is computed as in \cref{aniso}, with coefficients $k=\num{1e-6} \textrm{ m}^2$, $\beta = 100$ and $\alpha \in \left[-\pi/4,\pi/4\right]$. We assume a \textit{bottom} TL position, where the outer boundary of the TL overlaps the outer contour of the porous block (see \cref{figure_5_1}) and $d_TL = \num{1e-3}$ m. We again assume that the spatial variation of the porosity and the coefficients of the inverse of the permeability tensor within TL are given by \cref{eq:3}, with $d_{\epsilon} = d_K = 0$ and $\theta_{\epsilon}=\theta_K = 8000$. \par
The aim of this \textit{showcase} test is to compare the solution of the presented ODA solver, in terms of velocity and pressure fields, as well as fluxes crossing the boundary of the porous obstacle, with the numerical solution of a TDA solver proposed in \cite{46}, which couples the Navier Stokes equations (for compressible fluids) in $\Omega_{ff}$ to Darcy flow equation in $\Omega_{pm}$, enforcing conservation of mass and momentum across the by interface and by applying the the Beavers and Joseph slip condition at the interface \cite{3}. A staggered-grid finite volume method is applied to discretize the Navier-Stokes equations and a MPFA finite volume method, for the discretization of the Darcy equation. The MPFA scheme is suitable to simulating anisotropic problems in porous media and does not require the computational mesh to be $\mathbf{K}$-orthogonal to the principal anisotropy directions (i.e., no specific mesh alignment along the principal direction of the permeability tensor is required). The numerical TDA-MPFA procedure is implemented in the open-source software $\textrm{DuMu}^\textrm{x}$. More details can be found in \cite{46}. \par
The computational domain in the case of $Re\ll 1$ is $\Omega = \left[0, 0.75\right] \times \left[0,0.25\right] \textrm{ m}^2$, the porous obstacle $\Omega_{pm} = \left[0.25, 0.5\right] \times \left[0,0.2\right] \textrm{ m}^2$ with pressure drop $\Delta \Psi=\num {1e-6} \textrm{ m}^2/\textrm{s}^2$ (see \cref{figure_5_1}). 
The TDA-MPFA solver uses a structured grid with a uniform mesh size equal to $\num{5e-4}$ m. The adopted mesh sizes of the ODA solver runs are $\delta_0=\num{5e-3}$ m and $\delta_{TL}=\num{1e-4}$ m, $N = 235859$ and $N_T = 468685$, $\Delta t = 4$ s and the $CFL_{max}$ ranges from 1.13 ($\alpha = -\pi / 6 $) to 1.19 ($\alpha = \pi /4$).\par
\begin{figure} [h!]
	\centering
	\includegraphics[width=1\textwidth]{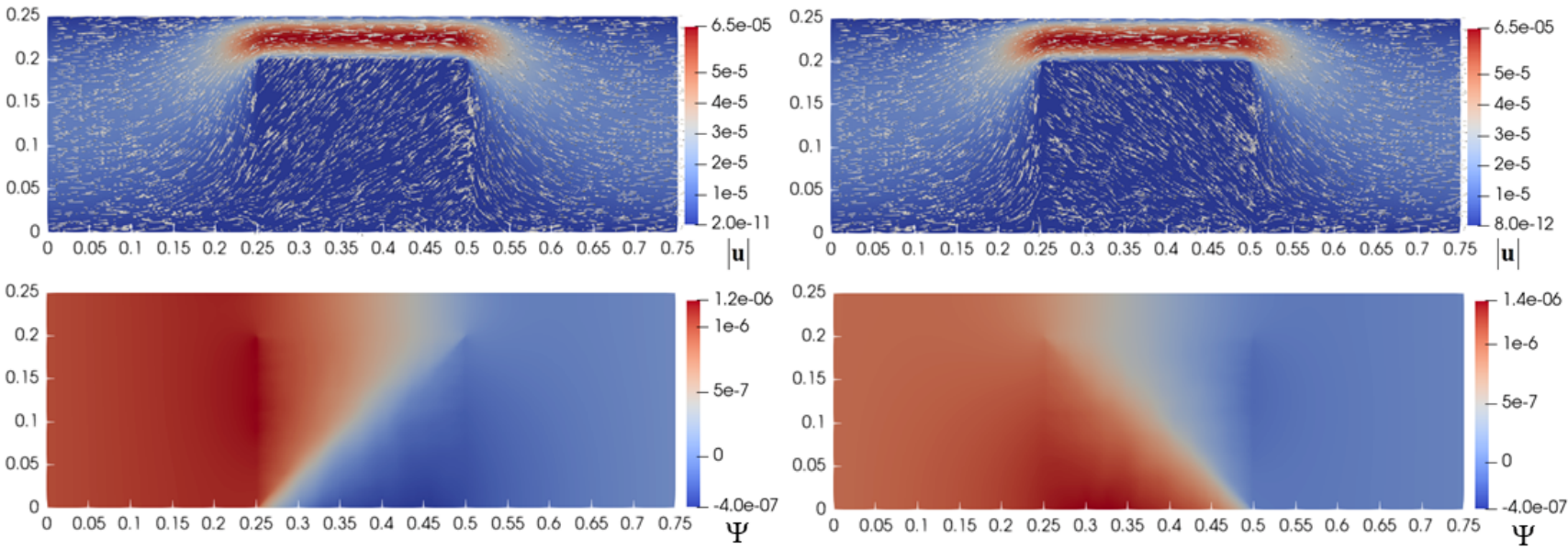}
	\caption{Test 5. Velocity and kinematic pressure fields computed by the ODA solver for the case of $Re \ll 1$. Left column: $\alpha=-\pi/4$, right column: $\alpha=\pi/4$. Top row: $\mathbf{u}$, bottom row: $\Psi$ (velocity in m/s, kinematic pressure in $\textrm{m}^2/\textrm{s}^2$)}
	\label{figure_5_2}
\end{figure}

 
In \cref{figure_5_2} we show the velocity and pressure fields provided by the presented ODA for $\alpha = -\pi/4$ and $\alpha = \pi/4$, in the case of $Re \ll 1$. Overall good agreement is observed with the results provided by the TDA-MPFA scheme in \cite{46} (see Fig. 5 of the referred paper). Due to the obstacle, a channelized flow is established above it, where the highest velocity values are observed, while the velocity in the porous block is approximately 2 - 4 orders of magnitude smaller. The effect of anisotropy is clearly visible within $\Omega_{pm}$, where the flow follows the principal direction of the permeability tensor, exiting ($\alpha=-\pi/4$) or entering ($\alpha=\pi/4$) at the top side. The recirculation zones simulated by the TDA-MPFA scheme within $\Omega_{pm}$ close to the bottom right corner ($\alpha = -\pi/4$) and left corner ($\alpha=\pi/4$) are slightly shifted outside $\Omega_{pm}$ in the presented ODA solutions. This could be caused by the different velocity distribution inside the TL. \par 		
For anisotropic problems, the TPFA scheme requires the computational grid to satisfy the $\mathbf{K}$-orthogonality (see Fig. 6 in \cite{46}). Otherwise the anisotropy is not correctly captured and the fluid flows almost horizontally in the porous block. \par
\begin{figure} [h!]
	\centering
	\includegraphics[width=0.7\textwidth]{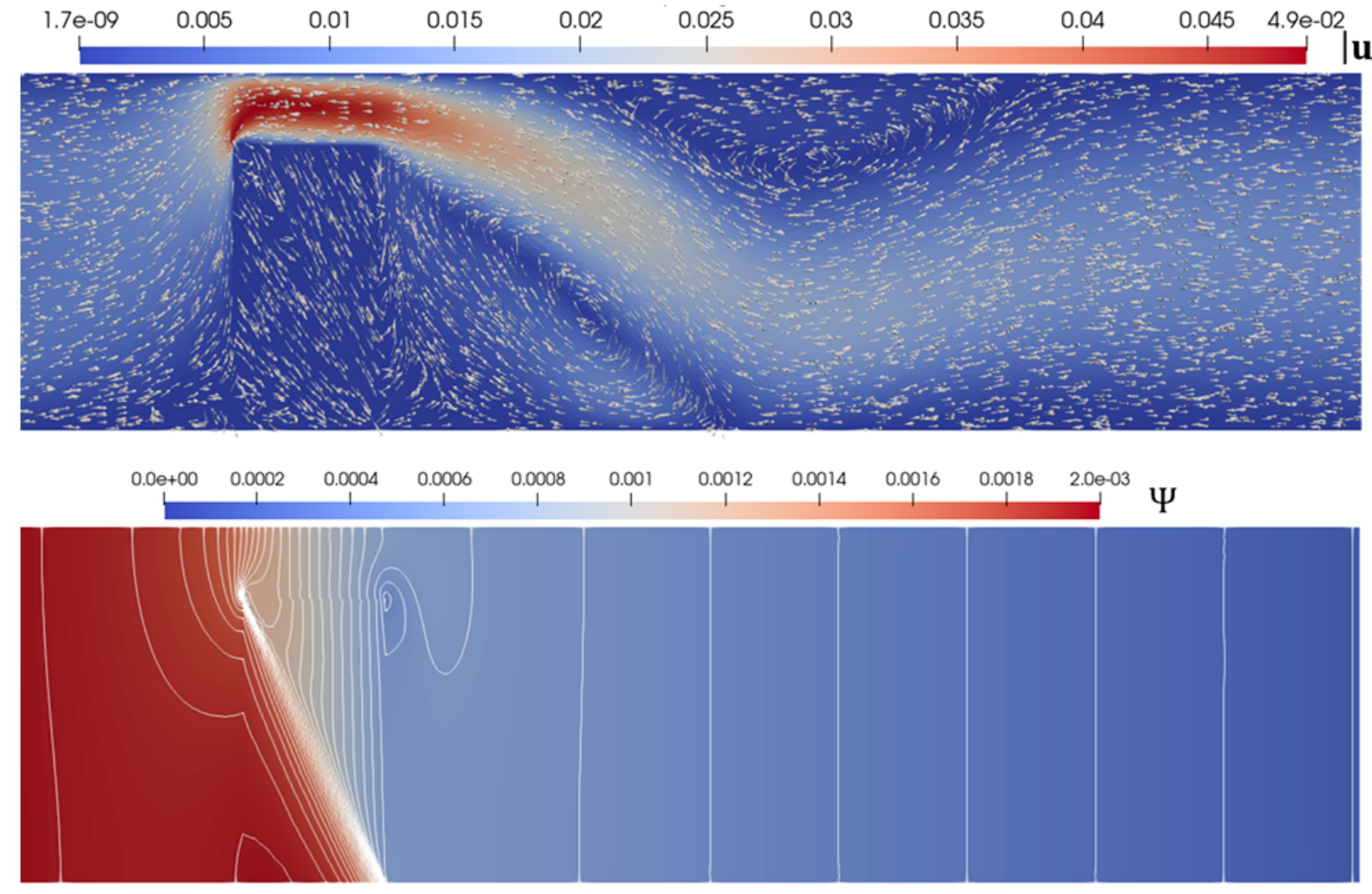}
	\caption{Test 5. Velocity and kinematic pressure fields computed by the ODA solver, case $Re \simeq 160$, $\alpha=\pi/4$. To improve visualization, the domain is scaled by a factor of 2 along the vertical direction (velocity in m/s, kinematic pressure in $\textrm{m}^2/\textrm{s}^2$)}
	\label{figure_5_4}
\end{figure}

 
In the case of $Re \simeq 160$, we adopt a similar setting to the previous one, with a longer domain along the $x$ direction. $\Omega = \left[0, 2.5\right] \times \left[0,0.25\right] \textrm{ m}^2$, the porous obstacle $\Omega_{pm} = \left[0.4, 0.6\right] \times \left[0,0.2\right] \textrm{ m}^2$ and pressure drop $\Delta \Psi=\num {2e-3} \textrm{ m}^2/\textrm{s}^2$. We investigate the case of $\alpha = \pi/4$. We use the same mesh sizes as before ($N= 344705$ and $N_T =680476$), the time step size is $\num{2.5e-2}$ s and $CFL_{max} \simeq 9.15$. We analyzed the case of $\alpha = \pi/4$. Again, the fluid is forced to flow mainly in the narrow channel over the porous obstacle, and vortex structures within $\Omega_{ff}$ are detected after approximately 20 s. The stationary solution is achieved after a longer time ($\sim$ 200 s, in \cref{figure_5_4}), compared to the case of $Re \ll 1$, with two stable countercurrent, large vortices downstream the porous block and a smaller one in front of the obstacle. Some discrepancies arise in the fluid region compared to the results provided in \cite{46} (see Fig. 8 in the referred paper). The reason could be the different assumption of compressible fluid made in the reference study. \par
\begin{figure}[h!]
	\centering
	\begin{subfigure}{\textwidth}
		\centering
		\includegraphics[width=\textwidth]{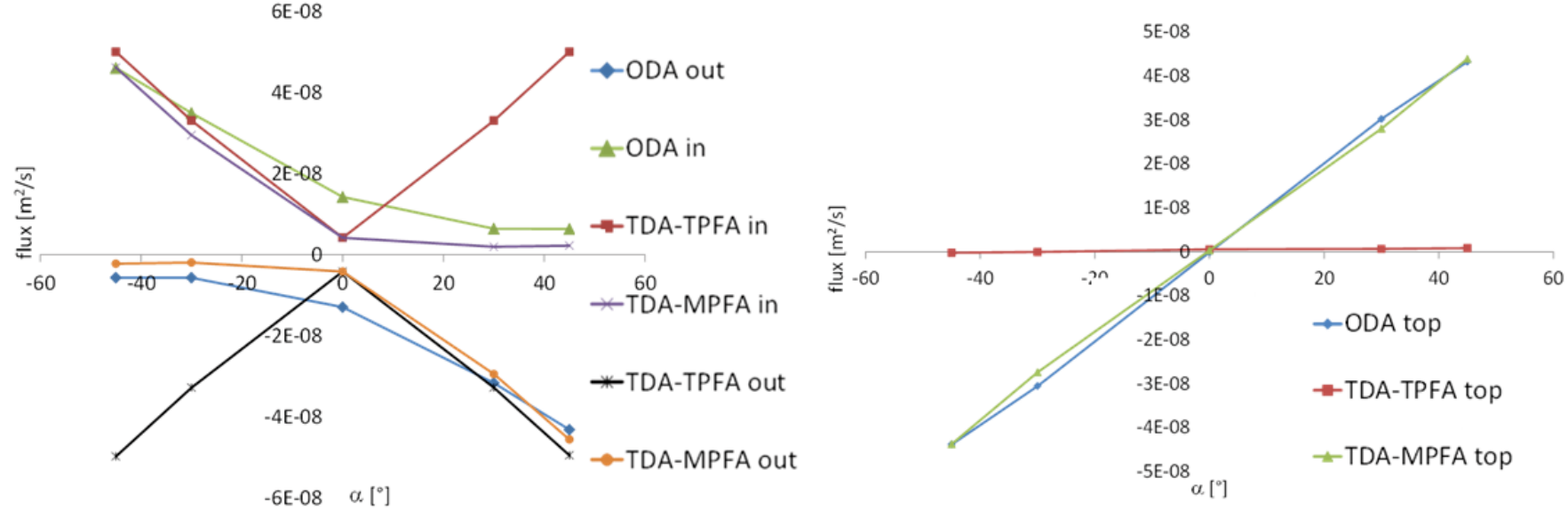}
		\caption{Case of $Re \ll 1$. Left) Fluxes on the upstream and downstream boundaries. Right) Fluxes at the upper boundary}
		\label{figure_5_3_a}%
	\end{subfigure}
	\hfill 
	\begin{subfigure}{0.8\textwidth} 
		\includegraphics[width=0.8\textwidth]{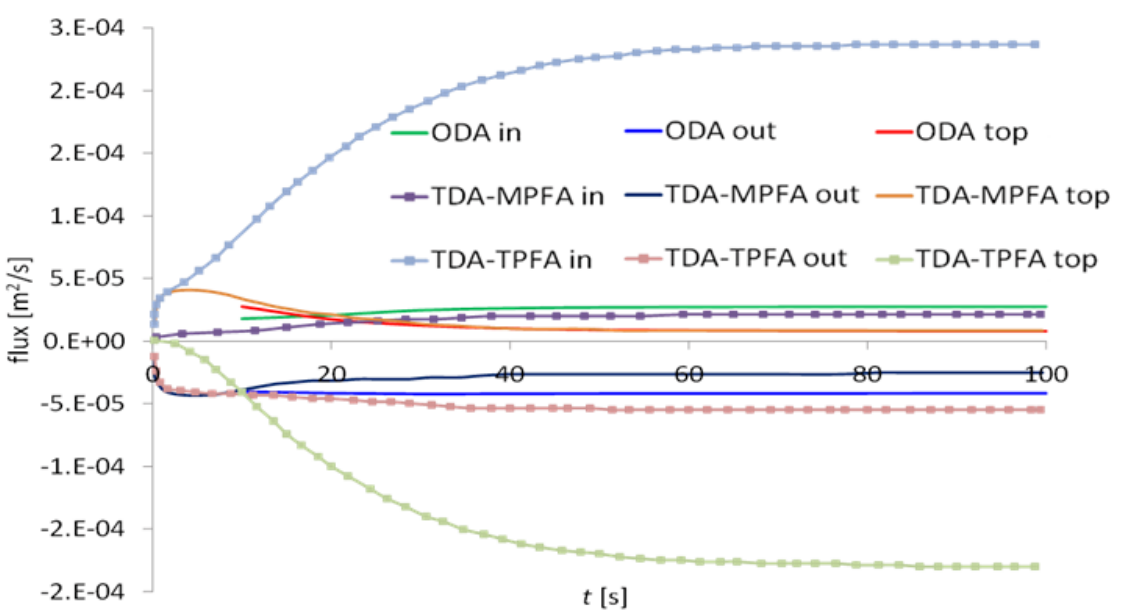}
		\caption{Case of $Re \simeq 160$. Time evolution of the fluxes on the upstream, downstream and upper boundaries}
		\label{figure_5_3_b}%
	\end{subfigure}
	\caption{Test 5. Computed fluxes crossing the boundary of $\Omega_{pm}$. Comparison between the ODA, TDA-MPFA and TDA-TPFA solvers (results of TDA-MPFA and TDA-TPFA from \cite{46})}
	\label{figure_5_3}
\end{figure}

  
It is interesting to compare the fluxes across the boundaries of the porous block computed by the three numerical solvers. The case of $Re \ll 1$ is plotted in \cref{figure_5_3_a}, where negative (positive) values are associated with fluxes leaving (incoming) the block. ODA and TDA-MPFA schemes provide similar results, and the discrepancies are due to the different treatment of the interface. The TDA-TPFA scheme computes correct results only for the $\mathbf{K}$-orthogonal case, i.e. when $\alpha = 0$. Since the off-diagonal coefficients of the tensor $\mathbf{K}$ are not considered in the TDA-TPFA scheme, the associated results are independent of the direction of rotation. In \cref{figure_5_3_b} we compare the time evolution of the fluxes crossing the boundary of the block for the case of $Re \simeq 160$. Again, ODA and TDA-MPFA solvers predict similar trends, with significant inflow crossing the top side of the obstacle, coming from the channel above the porous block, and a smaller amount of inflow through the downstream side of the block, due to the anisotropic effects. The results of the TDA-TPFA scheme do not match those of the two previous solvers, since again the anisotropic effects within $\Omega_{pm}$ are not properly captured.	

\section{Conclusions}  

A mesoscale ODA solver has been presented for the simulation of transfer processes between a free fluid and an anisotropic porous media. The governing equations are given by the Navier-Stokes-Brinkman equations together the continuity equation, assuming incompressible fluids. A fractional time step procedure is applied, by solving a prediction and two corrector steps within each time step. The numerical features and the associated advantages of the algorithm steps are presented and discussed. The numerical flux discretization strategy adopted in the corrector steps can be regarded as a Two-Point-Flux-Approximation (TPFA) scheme, but, unlike the standard TPFA schemes, the presented model correctly retains the anisotropy effects of the porous medium without the $\mathbf{K}$-orthogonality grid condition. This is proved by means of several tests. Very good agreement is obtained with both reference analytical and averaged pore-scale solutions. The proposed solver overcomes the restrictions of most other ODA solvers at the mesoscale that were recently presented in the literature, such as low Reynolds numbers, 1D flow or linearization of the convective inertial terms. In some of the numerical applications we discuss the discrepancies between the solutions provided by the present solver and a macroscopic-scale ODA algorithm, where a set of \textit{penalized} Navier-Stokes equations is solved. Compared with the reference solutions, the results of this \textit{penalized} algorithm show significant differences throughout the whole domain, not only close to the $\Omega_{pm} / \Omega_{ff}$ interface.

\section*{Funding}
This work was supported by the German Research Foundation (DFG) for supporting this work by funding SFB 1313, Project Number 327154368, Research Project A02, and by funding SimTech via Germany’s Excellence Strategy (EXC 2075 – 390740016).



\appendix
\renewcommand{\thesection}{A.\arabic{section}}
\setcounter{section}{0}
\renewcommand{\theequation}{\thesection.\arabic{equation}}

\section{Appendix A. Properties of the RT0 space functions} \label{App1}

\setcounter{equation}{0}

Any function $\mathbf{u}_e \in \mathfrak{R}_e$, where $\mathfrak{R}_e$ is the lowest-order Raviart-Thomas ($\mathbf{RT0}$) space function \cite{24}, is written as
\begin{equation} \label{eq:RT0}
	\mathbf{u}_e\left(\mathbf{x}\right) = \sum_{j=1}^3Q_j^e   \mathbf{w}_j^e  \quad \textrm{with} \quad \mathbf{w}_j^e=\frac{\left(\mathbf{x} - \mathbf{x}_j\right)}{2A_e} \quad \quad e=1,...,N_T,
\end{equation}
\noindent where $\mathbf{w}_j^e$ is the $j$-th space function of $\mathfrak{R}_e$, $\mathbf{x}_j$ is the coordinate vector of node $j$ in triangle $e$, opposite to side $j$, $A_e$ is the area of triangle $e$ and $Q_j^e$ is the flux crossing side $j$, positive outward. The properties of $\mathfrak{R}_e$ are 
\begin{subequations}  \label{eq:RT0_prop}
	\begin{gather}
		\nabla \cdot \mathbf{u}_e \textrm{ is constant inside each triangle } e  \\
		\mathbf{u}_e\cdot\mathbf{n}_j \textrm{ is constant for each side } j \in e,
	\end{gather}
\end{subequations}	
\noindent and $\mathbf{n}_j$ is the unit vector orthogonal to side $j$, pointing outward. According to \cref{eq:RT0}, the velocity components are piecewise linear inside each triangle $e$, and due to \cref{eq:RT0_prop}, $\mathbf{u}_e$ is piecewise constant within $e$ if $\sum_{j=1}^3Q_j^e=0$. If this condition is satisfied, $\nabla \cdot \mathbf{u}_e=0$  $\forall \mathbf{x}\in e$, $\forall e \in \Omega_T$, and, if the fluxes of two neighboring triangles are equal in value and opposite in sign along the common side, both local and global mass continuity are preserved. \par	

\section{Appendix B. Details of the numerical procedure applied for the CP1} \label{App2}

\setcounter{equation}{0}

We iteratively solve the two systems in \cref{eq:cs1_sys2} for the components of $\Delta \tilde{\mathbf{u}}_e$ unknowns, $\Delta \tilde{u}_e$ and $\Delta \tilde{v}_e$, respectively, as described in \cref{eq:A1}
\begin{subequations} \label{eq:A1}
	\begin{gather}
		\begin{cases} \label{eq:A1.1}
			i=0   \qquad  \textrm{! initialize the loop counter}\\			
			\begin{cases}
				\bar{\mathbf{M}}_{x,e} \begin{pmatrix}
					\Delta \tilde{u}_e^i \\ 0
				\end{pmatrix}A_e= \nu \sum_{j=1}^3 \ \frac{\Delta \tilde{u}_e^i-\Delta \tilde{u}_{ep}^i}{d_{e,ep}} l_j^e - \nu \sum_{j=1}^3 \frac{\Delta \breve{u}_e-\Delta \breve{u}_{ep}} {d_{e,ep}} l_j^e \\ \bar{\mathbf{M}}_{y,e} \begin{pmatrix}
					0 \\ \Delta \tilde{v}_{e}^i
				\end{pmatrix}A_e= \nu \sum_{j=1}^3 \ \frac{\Delta \tilde{v}_e^i-\Delta \tilde{v}_{ep}^i}{d_{e,ep}} l_j^e - \nu \sum_{j=1}^3 \frac{\Delta \breve{v}_e-\Delta \breve{v}_{ep}} {d_{e,ep}} l_j^e
			\end{cases}      \\
			\\
			\begin{cases} 
				\text{do while } (err > toll) \\
				\\
				\begin{cases}
					\bar{\mathbf{M}}_{x,e} \begin{pmatrix}
						\Delta \tilde{u}_e^{i+1} \\ \Delta \tilde{v}_e^i
					\end{pmatrix}A_e= \nu \sum_{j=1}^3 \ \frac{\Delta \tilde{u}_e^{i+1}-\Delta \tilde{u}_{ep}^{i+1}}{d_{e,ep}} l_j^e - \nu \sum_{j=1}^3 \frac{\Delta \breve{u}_e-\Delta \breve{u}_{ep}} {d_{e,ep}} l_j^e \\ \bar{\mathbf{M}}_{y,e} \begin{pmatrix}
						\Delta \tilde{u}_e^i  \\ \Delta \tilde{v}_e^{i+1}
					\end{pmatrix}A_e = \nu \sum_{j=1}^3 \ \frac{\Delta \tilde{v}_e^{i+1}-\Delta \tilde{v}_{ep}^{i+1}}{d_{e,ep}} l_j^e- \nu \sum_{j=1}^3 \frac{\Delta \breve{v}_e-\Delta \breve{v}_{ep}} {d_{e,ep}} l_j^e
				\end{cases}      \\
				\\
				i=i+1 \qquad \textrm{! update the loop counter}\\
				\text{end do}
			\end{cases}\\
		\end{cases},	\\
		\text{   with           } \bar{\mathbf{M}}_{x,e}=\begin{pmatrix} 
			\bar{M}_{1,1}^e\\\bar{M}_{1,2}^e
		\end{pmatrix}^T \qquad \bar{\mathbf{M}}_{y,e}=\begin{pmatrix}
			\bar{M}_{2,1}^e\\\bar{M}_{2,2}^e
		\end{pmatrix}^T \label{eq:A1.2},\\
		err=\min(err_x,err_y) \text{     with     } \begin{cases}
			err_x=\frac{\sqrt{\left(\| \Delta \tilde{u}_e^{i+1} -\Delta \tilde{u}_e^i   \|\right)}}{\sqrt{\|\Delta \tilde{u}_e^i\|}} \\ err_y=\frac{\sqrt{\left(\| \Delta \tilde{v}_e^{i+1} -\Delta \tilde{v}_e^i   \|\right)}}{\sqrt{\|\Delta \tilde{v}_e^i\|}}
		\end{cases} ,
	\end{gather}
\end{subequations}
\noindent where $i$ is the counter of the iterations in the iterative procedure and $1\text{d}-04\le toll \le1\text{d}-03$. \par 

\section{Appendix C. Details of the numerical procedure applied for the CP2} \label{App3}

\setcounter{equation}{0}
Starting from \cref{eq:A3.0},  
\begin{equation} \label{eq:A3.0}
	\left(\mathbf{\Xi}_j \nabla \eta\right) \cdot \mathbf{n}_j = \left( \mathbf{\Xi}_j \mathbf{n}_j \right) \cdot \nabla \eta,
\end{equation}
\noindent applying a co-normal decomposition, we obtain (see also \cref{fig_cs2})
\begin{equation} \label{eq:A3.1}
	\mathbf{\Xi}_j \mathbf{n}_j = \mathbf{d}_j \textrm{  with  } \mathbf{d}_j = \mathbf{d}_{j,n} + \mathbf{d}_{j,\tau} \textrm{  and  } \mathbf{d}_{j,\tau}=\mathbf{d}_{j,\tau_{n_1}} + \mathbf{d}_{j,\tau_{n_2}} ,
\end{equation}
\noindent where the vectors $\mathbf{d}_{j,n}$ and $\mathbf{d}_{j,\tau}$ are parallel to the directions $\mathbf{n}_j$ and $\boldsymbol{\tau}_j$, respectively, with $\boldsymbol{\tau}_j$ the unit vector tangential to side $l_j^e$. Let $\mathbf{n}_l$ and $\mathbf{n}_m$ (with \textit{l} = 1, 2 and \textit{m} = 3, 4), be the unit vector orthogonal to the sides shared by cells $e$ and $ep_l$, and by cells $ep$ and $ep_m$, respectively, pointing outwards from $e$ and $ep$, respectively (see \cref{fig_cs2}). According to the last relation in \cref{eq:A3.1}, vector $\mathbf{d}_{j,\tau}$ is decomposed along the $\mathbf{n}_1$ and $\mathbf{n}_2$ directions. \par
Starting from \cref{eq:A3.1}, we discretize the dot product $\mathbf{d}_j \cdot \nabla \eta$ in \cref{eq:Cs2_5} as  
\begin{subequations} \label{eq:A3.2}
	\begin{gather}
		\mathbf{d}_j \cdot \nabla \eta = \left(\mathbf{d}_{j,n} + \mathbf{d}_{j,\tau}\right) \cdot \nabla\eta \qquad \textrm{with} \\
		\mathbf{d}_{j,n} \cdot \nabla\eta = \frac{1}{2} \left(\frac{\eta_{ep} - \eta_e}{d_{e,ep}} - \frac{\eta_e - \eta_{ep}}{d_{e,ep}}\right) d_{j,n}  \qquad \textrm{and}\\
		\begin{gathered}
			\mathbf{d}_{j,\tau} \cdot \nabla\eta = \frac{1}{2} \left(\frac{\eta_{ep_1} - \eta_e}{d_{e,ep_1}} d_{\tau_{n_1}} \alpha_1 + \frac{\eta_{ep_2} - \eta_e}{d_{e,ep_2}} d_{\tau_{n_2}} \alpha_2\right) - \\
			\frac{1}{2} \left(\frac{\eta_{ep_3} - \eta_{ep}}{d_{ep,ep_3}} d_{\tau_{n_3}} \alpha_3 + \frac{\eta_{ep_4} - \eta_{ep}}{d_{ep,ep_4}} d_{\tau_{n_4}} \alpha_4 \right)
		\end{gathered},
	\end{gather}
\end{subequations}
\noindent where, with the help of \cref{fig_cs2}, $d_{j,n}=\mathbf{d}_j \cdot \mathbf{n}_j$, $d_{\tau_{n_{l(m)}}}=\mathbf{d}_{j,\tau} \cdot \mathbf{n}_{l(m)} $, $\eta_e$, $\eta_{ep}$, $\eta_{ep_{l(m)}}$ are the values of $\eta$ in the circumcenters of cells $e$, $ep$, and $ep_{l(m)}$, respectively, $d_{e,ep}$, $d_{e,ep_l}$ and $d_{e,ep_m}$ are the distances with a sign of the circumcenters of cells $e$ and $ep$, $e$ and $ep_l$, $ep$ and $ep_m$ respectively, and coefficient $\alpha_{l(m)}=1$ if $\mathbf{d}_{j,\tau} \cdot \mathbf{n}_{l(m)} > 0$ otherwise $\alpha_{l(m)}=-1$.\par 
According to \crefrange{eq:A3.0}{eq:A3.2}, \cref{eq:Cs2_5} becomes 
\begin{equation} \label{eq:A3.3}
	\begin{gathered}
		\sum_{j=1}^3 \left(\overline{Fl}_j^e \right) = \sum_{j=1}^3 \left(\frac{\eta_e - \eta_{ep}}{d_{e,ep}} d_{j,n} \right)l_j^e + \\
		\frac{1}{2}\sum_{j=1}^3 \left(\sum_{l=1,2} \frac{\eta_e - \eta_{ep_l}}{d_{e,ep_l}} d_{\tau_{n_l}} \alpha_l + 
		\frac{1}{2} \sum_{m=3,4} \frac{\eta_{ep_m} - \eta_{ep}}{d_{ep,ep_m}} d_{\tau_{n_m}} \alpha_m \right) l_j^e  
	\end{gathered},
\end{equation}
\noindent and \cref{eq:A3.3} form a system to be solved for the $\eta$ unknowns. \par

\bibliographystyle{elsarticle-num} 
\bibliography{biblio.bib}

\end{document}